\input amssym.tex
\magnification 1200

\newcount \version
\newcount \shortversion
\newcount \longversion
\shortversion 1\relax
\longversion 2\relax
% uncomment one of the two following lines to get a long or short version
%\version \shortversion\relax
\version \longversion\relax

% Use the following template:
% \ifnum \version=\shortversion
% SHORT
% \else
% LONG
% \fi

\font\smallrm=cmr8

\font\Bbb=msbm10
\def\BBB#1{\hbox{\Bbb#1}}
\font\smallBbb=msbm10 scaled 800
\def\smallbbb#1{{\hbox{\smallBbb#1}}}
\font\Frak=eufm10
\def\frak#1{{\hbox{\Frak#1}}}
\font\smallFrak=eufm10 scaled 800
\def\smallfrak#1{{\hbox{\smallFrak#1}}}

%Formulas
\def\vir{2.1}
\def\comm{2.2}
\def\degr{2.3}
\def\qas{2.4}
\def\omd{2.5}
\def\omdeg{2.6}
\def\Dab{2.7}
\def\Da{2.8}
\def\tens{2.9}
\def\glrel{2.10}
\def\loc{2.11}
\def\vfcom{2.12}
\def\vfmom{2.13}
\def\dxp{2.14}
\def\derd{2.15}
\def\dcom{2.16}
\def\trana{2.17}
\def\Yfm{3.1}
\def\dYfm{3.2}
\def\mtop{3.3}
\def\vaz{3.4}
\def\Yg{4.1}
\def\Yko{4.2}
\def\Yka{4.3}
\def\Yda{4.4}
\def\Ydo{4.5}
\def\Phida{4.6}
\def\vYdo{4.7}
\def\tomo{4.8}
\def\tomh{4.9}
\def\tomn{4.10}
\def\Cgg{4.11}
\def\Cdag{4.12}
\def\Cdako{4.13}
\def\Cdakb{4.14}
\def\Cdadb{4.15}
\def\Ckikj{4.16}
\def\Cgki{4.17}
\def\Cdog{4.18}
\def\Cdoko{4.19}
\def\Cdoka{4.20}
\def\Cdoda{4.21}
\def\Cdodo{4.22}

%Theorems
\def\VOX{2.1}
\def\maphi{2.2}
\def\trip{2.3}
\def\isom{2.4}
\def\locl{2.5}
\def\Pstruc{3.1}
\def\Mmap{3.2}
\def\sheafL{3.3}
\def\main{4.1}
\def\twom{4.2}
\def\chir{4.3}
\def\ratl{5.1}

\def\C{{\BBB C}}
\def\Z{{\BBB Z}}

\def\vac{{\bf 1}}

\def\L{{\cal L}}

\def\R{\Omega_{ch}}

\def\M{{\cal M}}
\def\G{{\cal G}}

\def\V{{\cal V}}
\def\O{{\cal O}}

\def\bom{{\overline{\Omega}}}
\def\g{\frak{g}}
\def\ss{\frak{s}}

\def\Map{\hbox{\rm Map}}
\def\Vect{\hbox{\rm Vect}}
\def\End{\hbox{\rm End}}
\def\Tr{\hbox{\rm Tr}}
\def\Id{\hbox{\rm Id}}
\def\res{\hbox{\rm res}}
\def\rank{\hbox{\rm rank} \;}
\def\dim{\hbox{\rm dim}}

\def\hX{{\widehat{X}}}
\def\hU{\widehat{U}}
\def\dg{{\g}}
\def\hg{\widehat{\g}}
\def\gl{{{\frak{gl}}_N}}
\def\hgl{{\widehat{{\frak{gl}}}_N}}
\def\bu{\overline{u}}
\def\d{\partial}
\def\vv{{\nu}}
\def\tx{\tilde{x}}
\def\td{\tilde{\d}}
\def\tu{\tilde{u}}
\def\tv{\tilde{v}}
\def\tw{\tilde{w}}
\def\tom{\tilde{\om}}
\def\tE{\tilde{E}}
\def\tI{\tilde{I}}
\def\hx{\hat{x}}
\def\hd{\hat{\d}}
\def\hu{\hat{u}}
\def\hv{\hat{v}}
\def\hE{\hat{E}}

\def\Hei{{\hbox{\Frak Hei}}}
\def\Vir{{\hbox{\Frak Vir}}}
\def\VHei{V_{\hbox{\smallFrak Hei}}}

\def\Vaff{V_{{\hbox{\smallFrak g}}}}
\def\VVir{V_{\hbox{\smallFrak Vir}}}
\def\VZN{V_{{\smallbbb Z}^N}}
\def\Maff{M_{{\hbox{\smallFrak g}}}}
\def\MVir{M_{\hbox{\smallFrak Vir}}}
\def\Laff{L_{{\hbox{\smallFrak g}}}}
\def\LVir{L_{\hbox{\smallFrak Vir}}}
\def\CHei{C_{\hbox{\smallFrak Hei}}}
\def\Cg{C_{{\hbox{\smallFrak g}}}}
\def\CVir{C_{\hbox{\smallFrak Vir}}}
\def\om{\omega}
\def\omvir{\om^{\hbox{\smallFrak Vir}}}
\def\omhei{\om^{\hbox{\smallFrak Hei}}}
\def\omg{\om^{\hbox{\smallFrak g}}}
\def\tomvir{\tilde{\om}^{\hbox{\smallFrak Vir}}}
\def\tomhei{\tilde{\om}^{\hbox{\smallFrak Hei}}}

\def\Vgl{V_{{{{\smallfrak{gl}}}_N}}}
\def\Mgl{M_{{{{\smallfrak{gl}}}_N}}}
\def\Lgl{L_{{{{\smallfrak{gl}}}_N}}}
\def\Cgl{C_{{{{\smallfrak{gl}}}_N}}}
\def\omgl{\om^{{{{\smallfrak{gl}}}_N}}}
\def\tomgl{\tilde{\om}^{{{{\smallfrak{gl}}}_N}}}

\def\tomg{\tilde{\om}^{\hbox{\smallFrak g}}}
\def\cvir{{c \, \hbox{\smallrm dim} \, \dg \over c + h^\vee}}

\def\dz{{z_1^{-1} \delta \left( {z_2 \over z_1} \right)}}
\def\ddz{{z_1^{-1} {\d \over \d z_2} \delta \left( {z_2 \over z_1} \right)}}

\def\dirlim{\lim_{\longrightarrow \atop U}}

\centerline
{\bf Modules for a sheaf of Lie algebras on loop manifolds.}

\centerline
{ {\bf Yuly Billig} \kern-3pt
\footnote{${}^1$}{School of Mathematics and Statistics, Carleton University,
1125 Colonel By Drive, Ottawa, K1S 5B6, Canada. E-mail: billig@math.carleton.ca}}
\footnote{}{2010 Mathematics Subject Classification: 17B66, 17B69.}

\

\

{\bf Abstract.}
 We consider a central extension of the sheaf of Lie algebras of maps from a manifold
$\C^* \times X$ into a finite-dimensional simple Lie algebra, together with the sheaf
of vector fields on $\C^* \times X$. Using vertex algebra methods we construct sheaves
of modules for this sheaf of Lie algebras. Our results extend the work of
Malikov-Schechtman-Vaintrob on the chiral de Rham complex.

\

\

{\bf 0. Introduction.}

 Two most interesting examples of infinite-dimensional Lie algebras, affine
Kac-Moody algebras and the Virasoro algebra, are associated with a circle as an
underlying geometric object. In this paper we are going to make a transition from 
the circle to more general manifolds. To construct an analogue of (untwisted) affine
Kac-Moody algebra in this case we start with the algebra of functions on a manifold
with values in a finite-dimensional simple Lie algebra $\dg$, then take its central 
extension and add the Lie algebra of vector fields on the manifold, acting as derivations.

Our goal is to develop a representation theory for this class of Lie algebras. Since we
would like to retain the features of the theory of the highest weight modules, we still
need the concept of positive/negative Fourier modes, and for this reason as the
underlying manifold we take $\hX = \C^* \times X$, where $X$ is a smooth
irreducible complex algebraic variety of dimension $N$. The punctured complex line
$\C^*$ here is a complex analogue of a circle. We choose to work with complex manifolds 
only as a matter of technical convenience, and one could just as well
consider  $\hX = S^1 \times X$, where $X$ is a real manifold.

When we look at functions on an algebraic manifold, taking the algebra of globally
defined functions may be inadequate (for example, in case of complex projective 
manifolds this algebra contains only constant functions). Instead, it is natural
to use the language of sheaves.

We begin by taking the sheaf $\Map (\hX, \dg)$ of functions on $\hX$ with values
in $\dg$, or, equivalently, functions on $X$ with values in the loop algebra
$\C [t, t^{-1}] \otimes \dg$. As the space of the central extension we take
the sheaf $\bom^1_{\hX}$ associated with the quotient of 1-forms by differentials
of functions $\Omega^1_{\hX} / d \O_{\hX}$. There is a 2-cocycle on
$\Map (\hX, \dg)$ with values in $\bom^1_{\hX}$ that naturally generalizes 
the central cocycle on loop Lie algebras.  Finally we take a semidirect product $\G$
of this central extension with the sheaf $\Vect(\hX)$ of vector fields on $\hX$.

To get a representation theory for this sheaf of Lie algebras we need to construct 
the sheaves of modules. This is done using vertex algebras. We introduce a sheaf of vertex 
algebras $\V$ on $X$ and show that the sheaf of Lie algebras $\G$ acts on $\V$.
We then construct sheaves of modules for $\V$, which automatically become modules
for $\G$.

If we take $\dg$ to be a trivial Lie algebra, $\dg = (0)$, we get representations for the
sheaf of vector fields on $\hX$. We would like to point out here a connection
with the important construction of the chiral de Rham complex. Chiral de Rham 
complex was introduced by Malikov-Schechtman-Vaintrob in [9]. It is a sheaf 
of vertex superalgebras on $X$ with a $\Z$-grading and a differential.
Malikov and Schechtman [8] show that the chiral de Rham complex admits the action
of the sheaf of Lie algebras $\C[t,t^{-1}] \otimes \Vect(X)$. Our result yield 
a stronger statement that in fact a larger sheaf $\Vect(\hX)$ acts on the chiral 
de Rham complex. In addition to this we note that the chiral differential is a 
homomorphism of $\Vect(\hX)$-modules.

 In classical differential geometry the Lie  algebra of vector fields acts on modules of 
tensor fields of a fixed type, and modules of differential forms appearing in the classical 
de Rham complex are a special case of this. Likewise the modules that we construct here 
for $\Vect(\hX)$ could be thought of as chiralizations of tensor modules and we get a 
wider  class of representations than those appearing in the chiral de Rham complex.

In case when $X$ is a torus, the representation theory of toroidal Lie algebras and Lie algebra of vector fields was developed in [10, 4, 6, 1, 2, 3]. Since
$(\C^*)^N$ can be covered with a single chart, there was no need to work with sheaves
of Lie algebras. In the toroidal case one gets strong results on irreducibility of the modules [2, 3].

 In the present paper we work with the algebraic varieties, however, this theory also
works if $X$ is taken to be an analytic or a $C^\infty$ manifold. The loop component
$\C^*$ should be still viewed in the algebraic setting with the ring of Laurent polynomials as the algebra of functions.

 The structure of this paper is as follows: in Section 1 we introduce the sheaf $\G$
of Lie algebras, generalizing the construction of affine Kac-Moody algebras, in Section 2
we construct a sheaf $\V$ of vertex algebras on $X$ and in Section 3 we
define the sheaves of the generalized Verma modules $\M$ and their quotients $\L$.
In Section 4 we prove that the sheaves $\V$, $\M$ and $\L$ are modules
for the sheaf $\G$ of Lie algebras. In the final section we consider a version of our construction in the setting of rational functions on a manifold.

\

{\bf Acknowledgements.}
 I am grateful to Fedor Malikov for the helpful discussions.
This work is supported in part with a grant from the Natural Sciences and Engineering Research Council of Canada.

\

{\bf 1. A sheaf of Lie algebras.}

Let $X$ be a smooth irreducible algebraic variety over $\C$ of dimension $N$.
Let $\hX = \C^* \times X$. Fix a finite-dimensional simple Lie algebra $\g$ with a symmetric invariant bilinear form $( \cdot | \cdot )$. We consider the sheaf 
$\Map (\hX, \g)$ 
of locally regular functions on $\hX$ with values in $\g$ (or, equivalently,
functions on $X$ with values in the loop algebra
$\C[t,t^{-1}] \otimes \g$). This becomes a sheaf of Lie algebras over $X$ with pointwise multiplication. 

Even though our base manifold is $\hX$, all
sheaves that we consider throughout the paper are over $X$, and for each open set $U \subset X$, we will consider functions defined over $\hU = \C^* \times U$.

 Next we are going to model in this setting the construction of affine Kac-Moody algebras from the loop Lie algebras. Let $\bom^1_{\hX}$ be the sheaf
associated with the presheaf $\Omega^1_{\hX}/ d\O_{\hX}$, where
$\O_{\hX}$ is the sheaf of rational functions on $\hX$ and
$\Omega^1_{\hX}$ is the sheaf of 1-forms on $\hX$.

 We define a central extension sheaf of Lie algebras 
$$ \Map (\hX, \g) \oplus \bom^1_{\hX},$$
where $\bom^1_{\hX}$ is central, while the new Lie bracket on $\Map (\hX, \g)$
is defined by
$$[f_1 \otimes g_1, f_2 \otimes g_2] = f_1 f_2 \otimes [g_1, g_2] +
(g_1 | g_2) \overline{ f_2 d f_1} ,$$
where $f_1, f_2$ are functions on an open set $\hU$, $g_1, g_2 \in \g$,
and  $\, {\vrule height 7pt depth -6.5pt width 0.4cm} \,$  
is the canonical projection 
$\Omega^1_{\hX} \rightarrow \Omega^1_{\hX}/ d\O_{\hX}$.

 The sheaf of vector fields $\Vect(\hX)$ acts on 
$\Map (\hX, \g) \oplus \bom^1_{\hX}$ and we can form a semidirect product sheaf
$$ \left( \Map (\hX, \g) \oplus \bom^1_{\hX} \right) \rtimes \Vect(\hX) .$$
Here the action of $\Vect_{\hX}$ on $\Map (\hX, \g)$ is the natural action of
vector fields on functions, while the action on $\bom^1_{\hX}$ is via Lie
derivative
$$ \eta (\overline{f_1 d f_2}) = \overline{\eta(f_1) df_2} + 
\overline{f_1 d\eta(f_2)}, $$
for $\eta \in \Vect_{\hX} (\hU)$, $f_1, f_2 \in \O_\hX (\hU)$.

 The variety $X$ admits a finite covering by open affine sets $\{ U_i \}$
where each $U_i$ has local (uniformizing) parameters
$x_1, \ldots, x_N \in \O_X (U_i)$ such that 
$\Omega^1_X (U_i)$ is a free $\O_X (U_i)$-module of rank $N$ with generators
$dx_1, \ldots, dx_N$ [11]. 

 An open covering $\{ U_i \}$ of $X$ yields an open covering $\{ \hU_i \}$
of $\hX$. We fix a local parameter $t$ on $\C^*$, and we identify
functions on $\C^*$ with $\C[t, t^{-1}]$, so that
$$\O_\hX (\hU_i) = \C[t, t^{-1}] \otimes \O_X (U_i) .$$

 Let $\G$ be a sheaf of Lie algebras on $X$ and let $\M$ be a sheaf of
vector spaces on $X$.

{\bf Definition.} A representation $(\rho, \M)$ of $\G$ is a sheaf morphism
$$\rho: \G \times \M \rightarrow \M ,$$
such that for every open set $U \subset X$, the map
$$\rho_U : \; \G (U) \times \M (U) \rightarrow \M (U)$$
is a representation of the Lie algebra $\G(U)$.

 The main goal of this paper is to construct representations
of the sheaf of Lie algebras 
$$\G = \left( \Map (\hX, \g) \oplus \bom^1_{\hX} \right) \rtimes \Vect(\hX) .$$

\

{\bf 2. A sheaf of vertex algebras.}

 We will construct representations using the vertex algebra techniques. Our main object will be a sheaf $\V$ of vertex algebras. Locally this sheaf will
be defined as the space of functions with values in a certain vertex algebra
$V$.

Let us recall the basic notions of the theory of the vertex operator algebras.
Here we are following [5] and [7].

{\bf Definition.} 
{
A vertex algebra is a vector space $V$ with a distinguished vector
$\vac$ (vacuum vector) in $V$, 
an operator $D$ (infinitesimal translation) on the space $V$, 
and a linear map $Y$ (state-field correspondence)
$$\eqalign{
Y(\cdot,z): \quad V &\rightarrow (\End V)[[z,z^{-1}]], \cr
a &\mapsto Y(a,z) = \sum\limits_{n\in\Z} a_{(n)} z^{-n-1} 
\quad (\hbox{where \ } a_{(n)} \in \End V), \cr} $$
such that the following axioms hold:

\noindent
(V1) For any $a,b\in V, \quad a_{(n)} b = 0 $ for $n$ sufficiently large;

\noindent
(V2) $[D, Y(a,z)] = Y(D(a), z) = {d \over dz} Y(a,z)$ for any $a \in V$;

\noindent
(V3) $Y(\vac, z) = \Id_V$;

\noindent
(V4) $Y(a,z) \vac \in V [[z]]$ and $Y(a,z)\vac |_{z=0} = a$ for any $a \in V$
\ (self-replication);

\noindent
(V5) For any $a, b \in V$, the fields $Y(a,z)$ and $Y(b,z)$ are mutually local, that is, 
$$ (z-w)^n \left[ Y(a,z), Y(b,w) \right] = 0, \quad \hbox{ for \ } 
 n  \hbox{\ sufficiently large} .$$

A vertex algebra $V$ is called a vertex operator algebra (VOA) if, in addition, 
$V$ contains a vector $\omega$ (Virasoro element) such that

\noindent
(V6) The components $L(n) = \omega_{(n+1)}$ of the field
$$ Y(\omega,z) = \sum\limits_{n\in\Z} \omega_{(n)} z^{-n-1} 
= \sum\limits_{n\in\Z} L(n) z^{-n-2} $$
satisfy the Virasoro algebra relations:
$$  [ L(n) , L(m) ] = (n-m) L(n+m) + \delta_{n,-m} {n^3 - n \over 12} 
(\rank V) \Id, \quad \hbox{where \ } \rank V \in \C;  
\eqno{(\vir)}$$

\noindent
(V7) $D = L(-1)$;

\noindent
(V8) $V$ is graded by the eigenvalues of $L(0)$:
$V = \mathop\oplus\limits_{n\in \Z} V_n$ with $L(0) \big|_{V_n} = n \Id$.
}

This completes the definition of a VOA.

As a consequence of the axioms of the vertex algebra  
we have the  following important commutator formula:
$$\left[ Y(a,z_1), Y(b,z_2) \right] =
\sum_{n \geq 0} {1 \over n!}  Y(a_{(n)} b, z_2)
\left[ z_1^{-1} \left( {\d \over \d z_2} \right)^n
\delta \left( {z_2 \over z_1} \right) \right] . \eqno{(\comm)}$$
 As usual, the delta function is
$$ \delta(z) = \sum_{n\in\Z} z^n .$$
By (V1), the sum in the right hand side of the commutator formula
is actually finite.

All the vertex operator algebras 
that appear in this paper have the gradings by 
non-negative integers: $V = \mathop\oplus\limits_{n=0}^\infty V_n$.
In this case the sum in the right hand side of the commutator 
formula (\comm) runs from $n=0$ to $n = \deg(a) + \deg(b) -1$,
because 
$$\deg(a_{(n)} b) = \deg(a) + \deg(b) -n - 1 ,\eqno{(\degr)}$$ 
and the elements of negative degree vanish.

% The commutator formula (\comm) may be written as the commutator
% relations between the components of the vertex operators:
% $$ [a_{(n)}, b_{(m)} ] =
% \sum\limits_{j\geq 0} \pmatrix{ n \cr j \cr}
% (a_{(j)} b)_{(n+m-j)} . \eqno{(\Borc)} $$

Another consequence of the axioms of a vertex algebra is the Borcherds'
identity:
$$
(a_{(k)} b)_{(n)} c 
= \sum\limits_{j\geq 0} 
(-1)^{k+j+1} \pmatrix{k \cr j \cr}
b_{(n+k-j)} a_{(j)} c +
\sum\limits_{j\geq 0} 
(-1)^j \pmatrix{k \cr j \cr}
a_{(k-j)} b_{(n+j)} c, \;\;\;  k,n \in \Z. \eqno{(\qas)} $$ 
%
% $$\sum\limits_{j\geq 0} \pmatrix{m \cr j \cr}
% (a_{(k+j)} b)_{(m+n-j)} c $$
% $$= \sum\limits_{j\geq 0} 
% (-1)^{k+j+1} \pmatrix{k \cr j \cr}
% b_{(n+k-j)} a_{(m+j)} c +
% \sum\limits_{j\geq 0} 
% (-1)^j \pmatrix{k \cr j \cr}
% a_{(m+k-j)} b_{(n+j)} c, \quad \quad  k,m,n \in \Z. \eqno{(\qas)} $$ 

%We will particularly need its special case when $m=0$ and $k = -1$:
%$$ (a_{(-1)} b)_{(n)} c = 
%\sum\limits_{j\geq 0} b_{(n-j-1)} a_{(j)} c 
%+ \sum\limits_{j\geq 0} 
%a_{(-1-j)} b_{(n+j)} c, \quad \quad  k \in \Z. \eqno{(\qass)} $$ 

Let us list some other consequences of the axioms of a vertex algebra that we
will be using in the paper. It follows from (V7) and (V8) that 
$$ \omega_{(0)} a = D(a) \eqno{(\omd)} $$
and
$$ \omega_{(1)} a = \deg (a) a, \quad \quad \hbox{\rm for \ } a
\hbox{\rm \  homogeneous}. \eqno{(\omdeg)} $$
The map $D$ is a derivation of the $n$-th product:
$$ D(a_{(n)} b) = (Da)_{(n)} b + a_{(n)} Db . \eqno{(\Dab)}$$
It could be easily derived from (V2) that
$$\left( Da \right)_{(n)} = -n a_{(n-1)}. \eqno{(\Da)}$$
% and thus
% $$ a_{(-1-k)} = {1\over k!} (D^k (a))_{(-1)} ,  \quad k \geq 0. 
% \eqno{(\Dka)}$$

\

 The vertex algebra $V$ that we need for the construction of the sheaf $\V$,
is the tensor product of four well-known vertex algebras:
$$ V =\VHei \otimes \Vgl \otimes \Vaff \otimes \VVir . \eqno{(\tens)} $$
Let us describe each tensor factor.

 Consider a Heisenberg Lie algebra $\Hei$ with the basis
$\{u_p (j), v_p (j), \CHei |^{p=1,\ldots, N}_{j\in\Z} \}$
and the Lie bracket
$$ [u_p (m), v_q (s)] = m \delta_{pq} \delta_{m, -s} \CHei ,$$
$$ [u_p (m), u_q (s)] = [v_p (m), v_q (s)] = 0, \;\;
p,q = 1, \ldots, N, \; m,s \in\Z, $$
and the element $\CHei$ being central.
 The vertex algebra 
$$\VHei = \C [u_p(-j), v_{p}(-j) |^{p=1,\ldots,N}_{j=1,2,\ldots} ]$$
is a Fock space module for this Heisenberg algebra, in which $\CHei$ acts as an identity operator and the raising operators $u_p(j), v_p(j)$ with $j \geq 1$
annihilate the highest weight vector $\vac$. 
The generating fields for this vertex algebra are
$$u_p (z) = Y(u_p (-1) \vac, z) =
\sum_{j\in\Z}  u_p(j) z^{-j-1}, $$
$$v_p (z) = Y(v_p (-1) \vac, z) =
\sum_{j\in\Z} v_p(j) z^{-j-1}, \; \; p=1,\ldots,N ,$$
with $u_p(0)$ and $v_p(0)$ acting on $\VHei$ trivially.

% In any vertex algebra the commutator relations between the fields can be
% expressed through $n$-th products using the commutator formula:
% $$[Y(a,z_1), Y(b,z_2)] = \sum_{n \geq 0} {1\over n!} Y(a_{(n)} b, z_2)
% z_1^{-1} \left( {\d \over \d z_2} \right)^n \delta \left(
% {z_2 \over z_1} \right), \eqno{(\comm)}$$
% where $\delta(z)$ is the formal delta function,
% $$\delta(z) = \sum_{j\in\Z} z^j .$$

% The non-trivial relation between the generating fields in the Heisenberg vertex algebra 
% $$[v_p (z_1), u_q(z_2)] = \delta_{pq} \ddz $$
% can be thus expressed as
% $$(v_p(-1)\vac)_{(0)} u_q(-1) \vac = 0, \quad$$
% $$(v_p(-1)\vac)_{(1)} u_q(-1) \vac = \delta_{pq}\vac, \quad$$
% $$(v_p(-1)\vac)_{(n)} u_q(-1) \vac = 0 \quad {\rm for \ } n \geq 2 .$$

 The Virasoro element in $\VHei$ is
$$\omhei = \sum_{p=1}^N v_p (-1) u_p (-1) \vac$$
and $\rank (\VHei) = 2N$.

 Consider next the affine Lie algebra
$$\hgl = \C[t,t^{-1}] \otimes \gl \oplus \C \Cgl$$
with the Lie bracket
$$[t^m \otimes A, t^s \otimes B ] =
t^{m+s} \otimes [A, B] + m \delta_{m,-s} \Tr (AB) \Cgl,$$
where $A, B \in \gl(\C)$ and $\Cgl$ is central.

 The second tensor factor $\Vgl$ in (\tens) is the universal enveloping vertex algebra for $\hgl$ at level $1$. It is a highest weight module for $\hgl$, with
the highest weight vector being annihilated by the subalgebra
$\C[t] \otimes \gl$ and $\Cgl$ acting as the identity operator.
As a vector space, it is realized as
$$\Vgl =  U( t^{-1} \C [t^{-1}]\otimes \gl) \otimes \vac. $$
The generating fields of this vertex algebra are
$$E_{ab} (z) = Y(E_{ab} (-1) \vac ,z) =
\sum_{j \in \Z} E_{ab}(j) z^{-j-1} ,\; \; a,b = 1, \ldots, N,$$
where $E_{ab}$ is a matrix with entry $1$ in position $(a,b)$ and zeros elsewhere, and $E_{ab}(j) = t^j \otimes E_{ab}$. It follows from this
formula that $(E_{ab} (-1) \vac)_{(n)} = E_{ab} (n)$. 

The commutator relations between the generating fields are encoded in
$n$-th products:
$$E_{ab}(0) E_{cd}(-1) \vac = \delta_{bc} E_{ad}(-1) \vac 
- \delta_{ad} E_{cb}(-1) \vac ,$$
$$E_{ab}(1) E_{cd}(-1) \vac = \delta_{ad} \delta_{bc} \vac ,$$
$$E_{ab}(n) E_{cd}(-1) \vac = 0 \quad \hbox{\rm for \ } n \geq 2 . 
\eqno{(\glrel)}$$ 

 We consider the following (non-standard) Virasoro element in $\Vgl$:
$$\omgl =  {1\over 2(N+1)} \left(  I(-1) I(-1) \vac
+ \sum_{a,b=1}^N E_{ab} (-1) E_{ba}(-1) \vac \right) + {1\over 2} I(-2) \vac ,$$
where $I$ is the identity matrix. The rank of  $\Vgl$ is $-2N$ (see [2] for details).
 
 The third tensor factor is the universal enveloping vertex algebra for the
affine Kac-Moody algebra
$$\hg = \C[t,t^{-1}] \otimes \dg \oplus \C \Cg$$
at level $c$. As a vector space $\Vaff = U(t^{-1} \C[t^{-1}] \otimes \dg) \otimes \vac$
with $\C[t] \otimes \dg$ annihilating the vacuum vector $\vac$.

For $g_1, g_2 \in \dg$, the $n$-th products in this case are:
$$g_1 (0) g_2(-1) \vac = [g_1, g_2] (-1) \vac, \;
g_1 (1) g_2(-1) \vac = c (g_1 | g_2) \vac, \;
g_1 (n) g_2(-1) \vac = 0 \quad \hbox{\rm for \ } n \geq 2 .$$
When the level $c$ is non-critical, $c \neq - h^\vee$, where $h^\vee$ is the dual
Coxeter number of $\dg$, the vertex algebra $\Vaff$ has a Virasoro element $\omg$
and its rank is 
$$\rank(\Vaff) =  {c \, \dim \dg \over c + h^\vee} .$$

 To define the last tensor factor, consider the Virasoro Lie algebra $\Vir$
with the basis $\{ L(j), \CVir | j\in\Z \}$ and Lie bracket
$$[L(m), L(s)] = (m-s) L(m+s) + {m^3 - m \over 12} \delta_{m,-s} \CVir, 
\;\; m,s \in\Z,$$
and $\CVir$ being a central element.

 The vertex algebra $\VVir$ is the universal enveloping vertex algebra for
the Virasoro Lie algebra where the central element $\CVir$ acts as a scalar $- \cvir$. 
It is a highest weight module for the Virasoro algebra in which 
the operators $L(j)$ with $j \geq -1$ annihilate the highest weight vector 
$\vac$.
As a space, it is realized as 
$$\VVir = U(\Vir^{(-)}) \otimes \vac,$$
where $\Vir^{(-)}$ is the subalgebra spanned by $L(j)$ with $j \leq -2$.

 The generator of this vertex algebra is $\omvir = L(-2) \vac$ and the generating field of this vertex  algebra is
$$ Y(\omvir,z) = \sum_{j\in\Z} L(j) z^{-j-2} .$$

 The Virasoro element of $V$ is the sum of the Virasoro elements of its tensor factors:
$$\omega = \omhei + \omgl + \omg + \omvir .$$
The rank of the Virasoro tensor factor was chosen in a way to make the total rank of $V$
to be $0$.

 The $n$-th products for the rank $0$ Virasoro element are:
$$\omega_{(0)} \omega = D \omega, \quad
\omega_{(1)} \omega = 2 \omega, \quad
\omega_{(n)} \omega = 0 \quad \hbox{\rm for \ } n \geq 2 .$$

 We begin the construction of the sheaf $\V$ with its local description. 
We have fixed a covering of $X$ with 
open affine sets admitting local parameters. Let $U_i$ be one of these open sets
with local parameters $x_1, \ldots, x_N$. 

 We set 
$$\V(U_i) = V\otimes \O_X (U_i). \eqno{(\loc)}$$ 
The fields $v_p(z)$, $u_q(z)$, $E_{ab}(z)$ are the ``chiralizations'' of the
vector fields, 1-forms and $(1,1)$-tensors on $X$ respectively, and transform
under the changes of coordinates (see (\trana) below).

 Let us define the vertex algebra structure on the space 
$V\otimes \O_X (U_i)$. The vertex algebra $V$ is embedded as subalgebra 
$V \otimes 1$ in $V \otimes \O_X (U_i)$. The state-field correspondence map
$Y$ on the elements of $V$ is defined as above, with the only difference
that the action of $v_p (0)$ is now defined as
$$ v_p (0) = {\d \over \d x_p} ,$$
while $u_p(0)$ still acts as zero.
  
 We define the state-field correspondence map on $\vac \otimes  \O_X (U_i)$ as
$$Y(\vac \otimes f, z) = 
\sum_{s \in \Z_+^N} {1 \over s!} \bu(z)^s \otimes {\d^s f \over \d x^s},$$
where $\bu_p (z)$ is an antiderivative of $u_p(z)$:
$$\bu_p (z) = \sum_{j\in\Z \backslash \{ 0 \} } {1 \over j} u_p(-j) z^j , $$
and in general,
$$Y(\vv \otimes f, z) = : Y(\vv\otimes 1, z)Y(\vac \otimes f, z) : ,$$
for $f \in \O_X (U_i) $, $\vv\in V$. Here and throughout the paper
we use the multi-index notation, 
for $s = (s_1, \ldots, s_N) \in \Z_+^N$ we set
$s! = s_1 ! \ldots s_N !$, ${\d^s \over \d x^s} =
\left({\d \over \d x_1}\right)^{s_1} \ldots \left({\d \over \d x_N}\right)^{s_N}$,
$ \bu(z)^s =  \bu_1(z)^{s_1} \ldots  \bu_N(z)^{s_N}$, etc.

 Note that for $f,h \in  \O_X (U_i)$, 
$$Y(\vac \otimes fh, z) = Y(\vac \otimes f, z) Y(\vac \otimes h, z).$$

{\bf Proposition \VOX.} $V \otimes \O_X (U_i)$ is a vertex algebra.

{\bf Proof.} We are going to apply the Existence theorem ([5], Theorem 4.5).
The infinitesimal translation map $D$ on $V \otimes \O_X (U_i)$ is defined
in the following way:
$$D(\vv \otimes f) = D(\vv) \otimes f + \sum_{p=1}^N u_p(-1) \vv \otimes {\d f \over \d x_p},$$
where $\vv\in V$, $f\in \O_X (U_i)$.

One has to verify that the generating fields for $V \otimes \O_X (U_i)$ are mutually local.
The only non-trivial relation is between $v_a (z)$ and $Y(\vac \otimes f,z)$. It is easy to
check that
$$[v_a(z_1), \bu_b (z_2)] = \delta_{ab} z_1^{-1} \sum_{j\in\Z\backslash \{ 0 \} }
\left( {z_2 \over z_1} \right)^j ,$$
which implies 
$$[v_a(z_1), \bu (z_2)^s] = s_a \bu (z_2)^{s-\epsilon_a}  z_1^{-1} \sum_{j\in\Z\backslash \{ 0 \} }
\left( {z_2 \over z_1} \right)^j ,$$
while
$$[v_a(z_1), \Id \otimes f] = z_1^{-1} \Id \otimes {\d f \over \d x_a} .$$
In the above, $\epsilon_a$ is an element of $\Z^N$ with $1$ in position $a$ and zeros elsewhere.
Combining these, we get
$$[v_a(z_1), Y(\vac \otimes f,z_2)] $$
$$= \sum_{s \in \Z_+^N} {s_a \over s!} \bu(z)^{s-\epsilon_a}  \otimes {\d^s f \over \d x^s}
z_1^{-1} \sum_{j\in\Z\backslash \{ 0 \} } \left( {z_2 \over z_1} \right)^j 
 + \sum_{s \in \Z_+^N} {1 \over s!} \bu(z)^s \otimes z_1^{-1} {\d^s \over \d x^s}
{\d f \over \d x_a}$$ 
$$ = \sum_{s \in \Z_+^N} {1 \over s!} \bu(z)^s \otimes {\d^s \over \d x^s}
{\d f \over \d x_a}  \dz 
 = Y(\vac \otimes {\d f \over \d x_a}, z_2)  \dz,  \eqno{(\vfcom)}$$ 
which implies locality. Verification of other conditions of the Existence theorem is straightforward.

\

 As a corollary of (\vfcom) we obtain
$$[v_a(k), f_{(m)}] = \left({\d f \over \d x_a}\right)_{(m+k)}. \eqno{(\vfmom)}$$

Also note that
$$u_p(-1) \vac = Dx_p .  \eqno{(\dxp)} $$

 Clearly, for any open subset $U \subset U_i$ we have a natural restriction
homomorphism of vertex algebras
$$\res_{U_i, U}: \;\; V\otimes \O_X (U_i) \rightarrow V\otimes \O_X (U) .$$

 Over the intersection $U_i \cap U_j$ we have defined two vertex algebra structures
on $V \otimes \O_X (U_i \cap U_j)$. We are now going to construct the gluing isomorphism
between these structures. Let $x_1, \ldots, x_N$ be the local parameters on $U_i$ and
$\tx_1, \ldots, \tx_N$ be the local parameters on $U_j$. In order to emphasize
the fact that the fields 
$v_p(z)$, $u_q(z)$, etc.,  transform under the change of  the local parameters, we will denote them as $\tv_p(z)$, $\tu_q(z)$, etc., when working
in coordinates $\tx_1, \ldots, \tx_N$. 

Under the coordinate changes, the partial derivative operators
$\d_k = {\d \over \d x_k}$, $\td_s = {\d \over \d \tx_s}$ transform in 
the standard way:
$$\td_a = (\td_a x_p)\d_p, \quad \d_b = (\d_b \tx_s) \td_s .$$
Throughout this paper we use Einstein notations on summation over repeated indices.

The product of the jacobians of the coordinate changes is the identity:
$$ \delta_{ab} = \d_b x_a = (\d_b \tx_s)(\td_s x_a), \quad
\delta_{ab} = \td_a \tx_b = (\td_a x_p)(\d_p \tx_b) .$$
Further differentiating the last equality, we get
$$0 = \d_q \delta_{ab} = 
(\d_q \td_a x_p)(\d_p \tx_b) + (\td_a x_p)(\d_q \d_p \tx_b)
= (\d_q \td_a x_p)(\d_p \tx_b) + \td_a \d_q \tx_b . \eqno{(\derd)}$$

The operators $\d_r$ and $\td_p$ do not commute and their commutator
may be expressed as follows:
$$[\d_r, \td_p] = (\d_r \td_p x_q) \d_q = - (\td_p \d_r \tx_s) \td_s . \eqno{(\dcom)}$$

We now define the gluing isomorphism $\Phi_{ij} : V \otimes \O_X (U_j \cap U_i)
\rightarrow V \otimes \O_X (U_i \cap U_j)$, where in the first copy we use local parameters 
$\tx_1, \ldots, \tx_N$, while in the second we use $x_1, \ldots, x_N$. This map will be first 
defined on the generators of this vertex algebra and then extended as a vertex algebra homomorphism.
In order to simplify notations we will drop from now on the tensor product symbol, as well as the symbol $\vac$, and write $v_a(-1) f$ instead of $v_a(-1) \vac \otimes f$, etc.
$$\eqalign{
&\Phi_{ij} (f) =  f, \;\;\;\;  f \in \O_X(U_i \cap U_j),  \cr
&\Phi_{ij} (\tu_a (-1) \vac) = 
u_p (-1) \d_p \tx_a , \cr
&\Phi_{ij} (\tv_a(-1) \vac) = 
v_p (-1) \td_a x_p  + E_{sp}(-1) \d_s \td_a x_p ,  \cr
&\Phi_{ij} (\tE_{ab} (-1) \vac) =
E_{sp}(-1)  (\d_s \tx_a) (\td_b x_p)
+ u_s(-1) \td_b \d_s \tx_a \cr
& \hbox{\hskip 2.2cm}
= E_{sp}(-1) (\d_s \tx_a) (\td_b x_p) - (\td_b x_p)_{(-2)} (\d_p \tx_a), \cr
&\Phi_{ij} (g(-1) \vac) = g(-1) \vac , \;\;\;  g\in\dg,  \cr
&\Phi_{ij} (\tomvir) = \omvir . \cr} \eqno{(\trana)}$$

{\bf Lemma \maphi.} The map $\Phi=\Phi_{ij}$ extends to a homomorphism of vertex algebras.

{\bf Proof.} We need to show that the fields corresponding to the images of 
the generators satisfy the same relations as the original generating fields.
Taking into account the commutator formula (\comm), we need to show that
$$\Phi (a)_{(n)} \Phi (b) = \Phi (a_{(n)} b)$$
for any pair $a, b$ of the generators and all $n \geq 0$.

 The relations between $\Phi(v_p(-1) \vac)$ and $\Phi(u_q(-1) \vac)$
have been established in ([9], Theorem 3.7). Let us verify other relations.

Let us check that
$$\Phi (\tv_m(-1)\vac)_{(0)} \Phi(f) = \Phi(\tv_m(0) f ) .$$

In the proof of this lemma we will use extensively the Borcherds' identity (\qas):
% $$(a_{(-1)} b)_{(n)} = \sum_{k=0}^\infty b_{(n-k-1)} a_{(k)}
% + \sum_{k=-1}^{-\infty}  a_{(k)} b_{(n-k-1)} .$$
% 
% Then
$$\Phi (\tv_m(-1) \vac)_{(0)} \Phi(f) =
\left( v_p(-1) \td_m x_p
+ E_{sp} (-1) \d_s \td_m x_p \right)_{(0)} f$$
$$= \big( \td_m x_p \big)_{(-1)} v_p(0) f =
 (\td_m x_p) \d_p f = \td_m f
 = \Phi( \tv_m(0) f) .$$
The relations  
$$\Phi (\tv_m(-1) \vac)_{(n)} \Phi(f) = 0$$
for $n \geq 1$ follow from the degree considerations since the degree 
of the left hand side is $-n$.

Let us show now that
$$\Phi(\tv_a(-1) \vac)_{(n)} \Phi(\tE_{bc}(-1) \vac) = 0 \quad
\hbox{\rm for \ } n \geq 0.$$
We have
$$\Phi(\tv_a(-1) \vac)_{(0)} \Phi(\tE_{bc}(-1) \vac) $$
$$=\left( v_p(-1) \td_a x_p + E_{sp} (-1) \d_s \td_a x_p \right)_{(0)}
\left( E_{qr} (-1) (\d_q \tx_b ) (\td_c x_r )
+ u_q(-1) \td_c \d_q \tx_b \right)  $$
$$ = ( \td_a x_p)_{(-1)} v_p(0) E_{qr} (-1)  (\d_q \tx_b) (\td_c x_r)
+ (\d_s \td_a x_p)_{(-1)} E_{sp}(0) E_{qr}(-1) 
(\d_q \tx_b ) (\td_c x_r )$$
$$+ ( \td_a x_p )_{(-1)} v_p(0) u_q(-1) \td_c \d_q \tx_b 
+ ( \td_a x_p)_{(-2)} v_p(1) u_q(-1) \td_c \d_q \tx_b $$
$$+ (\d_s \td_a x_p)_{(-2)} E_{sp}(1) E_{qr}(-1) (\d_q \tx_b) (\td_c x_r)$$
$$= E_{qr}(-1) (\td_a x_p) \d_p ((\d_q \tx_b) (\td_c x_r)) 
+ E_{sr}(-1) ( \d_s \td_a  x_p) (\d_p \tx_b)  (\td_c x_r )$$
$$- E_{qp}(-1) ( \d_s \td_a x_p) (\d_q \tx_b) (\td_c x_s )$$
$$+ u_q(-1) (\td_a x_p) (\d_p \td_c \d_q \tx_b)
+ u_q(-1) ( \d_q \td_a x_p) (\td_c \d_p \tx_b)
+ u_q(-1) ( \d_q \d_s \td_a x_p) (\d_p \tx_b) (\td_c x_s)$$
$$= E_{qr} (-1) \td_a ( (\d_q \tx_b)(\td_c x_r))
+ E_{qr} (-1)  (\d_q \td_a x_p) (\d_p \tx_b) (\td_c x_r)
- E_{qr}(-1) (\td_c \td_a x_r) (\d_q \tx_b)$$
$$ + u_q (-1) \td_a \td_c \d_q \tx_b
 + u_q (-1) (\d_q \td_a x_p)(\td_c \d_p \tx_b)
 + u_q (-1) (\td_c \d_q \td_a x_p)(\d_p \tx_b) = 0.$$
In the above we used the relation (\derd).

$$\Phi(\tv_a(-1) \vac)_{(1)} \Phi(\tE_{bc}(-1) \vac)$$
$$= \left( v_p(-1) \td_a x_p  + E_{sp} (-1) \d_s \td_a x_p \right)_{(1)}
\left( E_{qr} (-1)(\d_q \tx_b) (\td_c x_r)
+ u_q(-1) \td_c \d_q \tx_b \right)  $$
$$= ( \td_a x_p)_{(-1)} v_p(1) u_q(-1) \td_c \d_q \tx_b 
+ (\d_s \td_a x_p)_{(-1)} E_{sp}(1) E_{qr}(-1)(\d_q \tx_b) (\td_c x_r)$$
$$ = (\td_a x_p)( \td_c \d_p \tx_b) + ( \d_s \td_a x_p)(\d_p \tx_b)(\td_c x_s)
 = \td_c \left( ( \td_a x_p)(\d_p \tx_b) \right) = 0 .$$

 The relations
$$\Phi(\tv_a(-1)\vac)_{(n)} \Phi(\tE_{bc}(-1)\vac) = 0 \quad
\hbox{\rm for \ } n \geq 2$$
follow trivially from the degree considerations (\degr).

 Let us now consider the analogues of (\glrel):
$$\Phi(\tE_{ab}(-1) \vac)_{(0)} \Phi(\tE_{cd}(-1) \vac)$$
$$= \left( E_{sp}(-1)(\d_s \tx_a ) (\td_b x_p)
+ u_s(-1) \td_b \d_s \tx_a \right)_{(0)}
\left( E_{qr}(-1)(\d_q \tx_c) (\td_d x_r)
+ u_q(-1) \td_d \d_q \tx_c \right)$$
$$ = \big((\d_s \tx_a) (\td_b x_p) \big)_{(-1)} E_{sp}(0) E_{qr}(-1)
(\d_q \tx_c)(\td_d x_r)$$
$$ + \big((\d_s \tx_a) (\td_b x_p) \big)_{(-2)}
E_{sp}(1) E_{qr}(-1)(\d_q \tx_c) (\td_d x_r)$$
$$= E_{sr}(-1)(\d_s \tx_a) (\td_b x_p) (\d_p \tx_c) (\td_d x_r)
- E_{qp}(-1)(\d_s \tx_a) (\td_b x_p) (\d_q \tx_c) (\td_d x_s)$$
$$+ u_j (-1) \d_j ((\d_s \tx_a) (\td_b x_p)) (\d_p \tx_c) (\td_d x_s)$$
$$ = \delta_{bc} E_{sr}(-1) (\d_s \tx_a) (\td_d x_r)
- \delta_{ad} E_{qp}(-1) (\d_q \tx_c) (\td_b x_p)
- \delta_{ad} u_j(-1) \td_b \d_j \tx_c
+ \delta_{bc} u_j (-1) \td_d \d_j \tx_a $$
$$ = \delta_{bc} \Phi(\tE_{ad}(-1) \vac) - \delta_{ad} \Phi(\tE_{cb}(-1) \vac) ,$$
and
$$\Phi(\tE_{ab}(-1) \vac)_{(1)} \Phi(\tE_{cd}(-1) \vac)
 = \big((\d_s \tx_a) (\td_b x_p) \big)_{(-1)} E_{sp}(1) E_{qr}(-1)
(\d_q \tx_c)(\td_d x_r) $$
$$=(\d_s \tx_a) (\td_b x_p) (\d_p \tx_c)(\td_d x_s)
= \delta_{ad} \delta_{bc} .$$

 Verification of the remaining relations is trivial.

\

{\bf Lemma \trip.} Over the triple intersection $U_i \cap U_j \cap U_k$ we have
$\Phi_{ij} \circ \Phi_{jk} = \Phi_{ik}$.

{\bf Proof.}
Let us denote $\hx_1, \ldots, \hx_N$ the local parameters on $U_k$. We need to verify
the equality $\Phi_{ij} \circ \Phi_{jk} = \Phi_{ik}$ on the generators of the vertex algebra
$V \otimes \O_X (U_k \cap U_j \cap U_i)$. For the functions $f \in  \O_X (U_k \cap U_j \cap U_i)$
this equality holds since the functions do not transform under the coordinate changes. Let us carry out 
the calculations for the generators of $V$.

$$\Phi_{ij} \left( \Phi_{jk}(\hu_a(-1) \vac) \right) = \Phi_{ij} ( \tu_p(-1) \td_p \hx_a)$$
$$= (u_s(-1) \d_s \tx_p)_{(-1)} \td_p \hx_a = u_s(-1) (\d_s \tx_p) (\td_p \hx_a)$$
$$= u_s(-1) \d_s \hx_a = \Phi_{ik}(\hu_a(-1) \vac) .$$

$$\Phi_{ij} \left( \Phi_{jk}(\hv_a(-1) \vac) \right) = \Phi_{ij} \left( \tv_p(-1) \hd_a \tx_p 
+ \tE_{sp}(-1) \td_s \hd_a \tx_p \right)$$
$$= \big( v_q(-1) \td_p x_q \big)_{(-1)} \hd_a \tx_p
+ \big(E_{rq}(-1) \d_r \td_p x_q \big)_{(-1)} \hd_a \tx_p$$
$$+ \big( E_{rq}(-1) (\d_r \tx_s) (\td_p x_q) \big)_{(-1)} \td_s \hd_a \tx_p
+ \left( u_r(-1) \td_p \d_r \tx_s \right)_{(-1)} \td_s \hd_a \tx_p$$
$$=v_q(-1) (\td_p x_q)(\hd_a \tx_p)
+ (\td_p x_q)_{(-2)} v_q(0) \hd_a \tx_p
+ E_{rq}(-1) (\d_r \td_p x_q)(\hd_a \tx_p)$$
$$+ E_{rq}(-1) (\d_r \tx_s) (\td_p x_q)(\td_s \hd_a \tx_p)
+ u_r(-1) (\td_p \d_r \tx_s)(\td_s \hd_a \tx_p)$$
$$=v_q(-1) \hd_a x_q 
+ u_r(-1) (\d_r \td_p x_q)(\d_q \hd_a \tx_p)
+ u_r(-1) (\td_p \d_r \tx_s) (\td_s \hd_a \tx_p)$$
$$+ E_{rq}(-1) (\d_r \td_p x_q)(\hd_a \tx_p)
+ E_{rq}(-1) (\td_p x_q)(\d_r \hd_a \tx_p)$$
$$= v_q(-1)\hd_a x_q + E_{rq}(-1) \d_r \hd_a x_q =  \Phi_{ik}(\hv_a(-1) \vac) .$$

In the above one can use (\dcom) to see that the terms with $u_r(-1)$ cancel.

$$\Phi_{ij} \left( \Phi_{jk}(\hE_{ab}(-1) \vac) \right) =
\Phi_{ij} \left( \tE_{sp}(-1)(\td_s \hx_a) (\hd_b \tx_p) + \tu_s(-1) \hd_b \td_s \hx_a \right)$$
$$=\left( E_{rq}(-1)(\d_r \tx_s)(\td_p x_q) + u_r(-1) \td_p \d_r \tx_s \right)_{(-1)}
(\td_s \hx_a)(\hd_b \tx_p)
+ \left( u_r(-1) \d_r \tx_s \right)_{(-1)} \hd_b \td_s \hx_a$$
$$= E_{rq}(-1)(\d_r \tx_s) (\td_p x_q)(\td_s \hx_a)(\hd_b \tx_p)
+ u_r(-1) (\hd_b \d_r \tx_s)( \td_s \hx_a)
+ u_r(-1) (\d_r \tx_s)(\hd_b  \td_s \hx_a)$$
$$= E_{rq}(-1)(\d_r \hx_a) (\hd_b x_q) 
+ u_r(-1) \hd_b \d_r \hx_a
=  \Phi_{ik}(\hE_{ab}(-1) \vac) .$$

\

{\bf Corollary \isom.} The map $\Phi_{ij}$ is an isomorphism.

{\bf Proof.} Setting $k=i$ in the previous Lemma, we get $\Phi_{ij}^{-1} =
\Phi_{ji}$.

\

 As a result of our construction, we get the following

{\bf Theorem \locl.} The local data (\loc) together with the gluing maps $\Phi_{ij}$
define a sheaf $\V$ of vertex algebras over $X$. 

\

{\bf 3. A sheaf of modules of chiral tensor fields.}

In this section we will construct sheaves of modules for the sheaf $\V$ of vertex
algebras.  

 First let us discuss the local situation. Let $U$ be an open set contained in $U_i$, 
and let $M = \mathop\oplus\limits_{n=0}^\infty M_n$ be a module for the vertex Lie algebra $\V(U)$.
  The module $M$ is a module for the Lie algebra
$$\ss = \Hei \oplus \hgl \oplus \hg \oplus \Vir .$$
This Lie algebra has a natural $\Z$-grading $\ss = \mathop\oplus\limits_{n\in\Z} \ss_n$ 
and a triangular decomposition $\ss = \ss_+ \oplus \ss_0 \oplus \ss_-$ associated with this grading.

 Using the standard methods of vertex algebras, we get the following

{\bf Proposition \Pstruc.} (i) $M_0$ is a module for  the commutative algebra $\O_X (U)$
with the action
$$ f m = f_{(-1)} m, \;\;\;  f \in \O_X (U), m \in M_0 .$$

(ii) $M_0$ is a module for the Lie algebra $\ss_0$.
 The actions of $\O_X (U)$ and $\ss_0$ on $M_0$ are compatible in the following 
way:
$$v_p(0) f m - f v_p(0) m = (\d_p f) m ,$$
while the remaining basis elements of $\ss_0$ commute with $\O_X (U)$.

(iii) For $f \in \O_X (U)$ introduce  the operator $T(f,z)$ on the space $U(\ss_-) M_0$:
$$T(f,z) y m = \sum_{k\in\Z_+^N} {1\over k!} \bu (z)^k y (\d^k f) m, $$
where $y \in U(\ss_-)$, $m \in M_0$.
 If $M$ is generated by $M_0$ as a $\V(U)$-module then 
$$Y(f,z) y m = T(f,z) y m .$$

(iv) If $M$ is generated by $M_0$ as a $\V(U)$-module then $M = U(\ss_-) M_0$.

{\bf Proof.} Part (i) follows from the relation
$$Y(f,z) Y(h,z) = Y(fh,z)$$
and the fact that $f_{(n)} m = 0$ for $n \geq 0$, $m\in M_0$.
Part (ii) is a consequence of the statement that $M$ is a graded $\ss$-module
and (\vfmom).
 
 We shall prove part (iii) by induction on the degree of $y$. For the basis
of induction, $\deg(y) = 0$, so that $y =1$, we need to show that
$$Y(f,z) m = \sum_{k\in\Z_+^N} {1\over k!} \bu (z)^k (\d^k f) m , \eqno{(\Yfm)}$$
for $m \in M_0$. It is clear that both sides involve only non-negative powers
of $z$. Let us reason by induction on $n$ the equality of terms up to $z^n$ 
in (\Yfm). The coefficients at $z^0$ in (\Yfm) coincide by the definition
of the action of $\O_X(U)$ on $M_0$. The equality of the coefficients at $z^{n+1}$ 
in (\Yfm) will follow from the equality of $z^n$ terms in
$$ {\d \over \d z} Y(f,z) m = {\d \over \d z} T(f,z) m. \eqno{(\dYfm)}$$
Note that
$$ {\d \over \d z} Y(f,z) m = \sum_{p=1}^N u_p(z) Y(\d_p f, z) m, $$
and also
$${\d \over \d z} T(f,z) m = \sum_{p=1}^N u_p(z) T(\d_p f, z) m. $$
Since in the above terms in $u_p(z)$ with the negative powers of $z$ act trivially
(note that $u_p(0)$ acts trivially on $M$ since $u_p(z) = {d \over d z} Y(x_p, z)$), 
the equality of $z^n$ terms in (\dYfm) follows from the induction assumption.

 Let us now complete the induction on the degree of $y \in U(\ss_-)$. Suppose
$y = y^\prime y^{\prime\prime}$, where $y^\prime \in \ss_-$,
$y^{\prime\prime} \in U(\ss_-)$. If $y^\prime$ is one of $u_a(-n)$, 
$E_{ab}(-n)$ or $L(-n)$, then both operators $Y(f,z)$ and $T(f,z)$ commute with
$y^\prime$ and we get
$$Y(f,z) y^\prime y^{\prime\prime} m =
y^\prime Y(f,z) y^{\prime\prime} m 
 = y^\prime  T(f,z) y^{\prime\prime} m 
= T(f,z) y^\prime y^{\prime\prime} m.$$
The only non-trivial case is $y^\prime =  v_a(-n)$. However it follows from (\vfcom)
that 
$$ [v_a(-n), Y(f,z) ] = z^{-n} Y(\d_a f,z),$$
and also
$$[v_a(-n), T(f,z) ] = z^{-n} T(\d_a f,z),$$
thus
$$Y(f,z) v_a(-n) y^{\prime\prime} m = v_a(-n) Y(f,z) y^{\prime\prime} m
- z^{-n} Y(\d_a f, z) y^{\prime\prime} m$$
$$ = v_a(-n) T(f,z) y^{\prime\prime} m - z^{-n} T(\d_a f, z) y^{\prime\prime} m
= T(f,z) v_a(-n) y^{\prime\prime} m .$$

 Part (iv) follows immediately from (iii).

\

{\bf Corollary \Mmap.} Let $M^\prime$, $M^{\prime\prime}$ be two $\V(U)$-modules,
$M^\prime = \mathop\oplus\limits_{n=0}^\infty M_n^\prime$,
$M^{\prime\prime} = \mathop\oplus\limits_{n=0}^\infty M_n^{\prime\prime}$,
that are generated by $M_0^\prime$ and $M_0^{\prime\prime}$ respectively. 
Let $\psi: M^\prime \rightarrow M^{\prime\prime}$ be a homomorphism
of $\ss$-modules preserving the grading, such that 
$\psi:  M_0^\prime \rightarrow M_0^{\prime\prime}$ is a homomorphism of
$\O_X (U)$-modules. Then $\psi$ is a homomorphism of $\V(U)$-modules.  

\

 The previous Proposition essentially tells us how the vertex algebra $\V(U)$ may act
on its modules. Let us now give an explicit construction.

 Let $W$ be a rational finite-dimensional simple $GL_N(\C)$-module and let
$\Mgl (W)$ be the generalized Verma module at level $1$ for $\hgl$, induced from
the $\gl$-module $W$. 
Let $\Maff(S)$ be the generalized Verma module for $\hg$ at level $c$,
induced from an irreducible $\dg$-module $S$, and $\MVir(h)$ be 
the Verma module for the Virasoro Lie algebra of rank $ -\cvir$ with the highest weight vector $v_h$ such that $L(0) v_h = h v_h$, $h \in \C$ .

  For an open set $U \subset U_i$, the space
 $$\M(U) =\VHei \otimes \Mgl (W) \otimes \Maff(S) \otimes \MVir(h) \otimes \O_X(U) $$
has a natural structure of a module for the vertex algebra $\V(U)$.

 In order to construct the sheaf of modules $\M$, we need to define the gluing
isomorphisms
$$\Psi_{ij} : \; \M(U_j \cap U_i) \rightarrow  \M(U_i \cap U_j) .$$
Both $\M(U_j \cap U_i)$ and $\M(U_i \cap U_j)$ are the modules for the vertex
algebra $\V(U_j \cap U_i)$, where the action on the second module is defined via the isomorphism  $\Phi_{ij}$ of vertex algebras. The isomorphism $\Psi_{ij}$ of modules that we need to construct must be an isomorphism of 
$\V(U_j \cap U_i)$-modules. The map $\Psi_{ij}$ will be constructed 
using Corollary \Mmap.
Note that both $\M^\prime = \M(U_j \cap U_i)$ and $\M^{\prime\prime} = 
\M(U_i \cap U_j)$ are $\Z$-graded and their top components are:
$$\M_0^\prime = W \otimes S \otimes v_h \otimes \O_X(U_j \cap U_i),$$
$$\M_0^{\prime\prime} = W \otimes S \otimes v_h \otimes \O_X(U_i \cap U_j).$$
We first construct the map $\Psi_{ij} : \;  \M_0^\prime \rightarrow \M_0^{\prime\prime}$, which is an isomorphism of $ \O_X(U_j \cap U_i)$-modules.
Note that the jacobian matrix $J^{(ij)} = J = (\d_r \tx_s) E_{rs} $ 
% with the entries  $J_{rs} = \d_r \tx_s$ 
is an element of $GL_N(\O_X(U_j \cap U_i))$. We set
$$\Psi_{ij} ( \tw \otimes s \otimes v_h \otimes f ) = (J^{(ij)} \tw) \otimes s \otimes
v_h \otimes f . \eqno{(\mtop)}$$
We claim that this is a homomorphism of $\ss_0$-modules.
Indeed, the action of $\dg$ on both $\M_0^\prime$ and $\M_0^{\prime\prime}$
is the natural action on $S$, while $L(0)$ acts as multiplication by $h$ on both spaces.
To see that   $\Psi_{ij} : \;  \M_0^\prime \rightarrow \M_0^{\prime\prime}$
is a homomorphism of $\gl$-modules, we need to check that
$$\Phi_{ij}(\tE_{ab}(0)) \Psi_{ij} (\tw \otimes s \otimes v_h \otimes f) =
 \Psi_{ij} (\tE_{ab}(0) \tw \otimes s \otimes v_h \otimes f) .$$

Applying (\trana), we see that the left hand side equals
$$E_{sp}(0) J \tw \otimes s \otimes v_h \otimes (\d_s \tx_a) (\td_b x_p) f,$$
and to compute the right hand side we use the connection between the action of the group $GL_N$ and its Lie algebra $\gl$ on $W$: 
$$J E_{ab} (0) \tw  \otimes s \otimes v_h \otimes f
= (J E_{ab} (0) J^{-1})  J \tw  \otimes s \otimes v_h \otimes f$$
$$= E_{sp}(0) J \tw \otimes s \otimes v_h \otimes (\d_s \tx_a) (\td_b x_p) f .$$

Since $u_a(0)$ acts on both $\M_0^\prime$ and $\M_0^{\prime\prime}$ trivially,
the last thing to check is the equality
$$\Phi_{ij}(\tv_a(0)) \Psi_{ij} (\tw \otimes s \otimes v_h \otimes f) =
 \Psi_{ij} (\tv_a(0) \tw \otimes s \otimes v_h \otimes f) . \eqno{(\vaz)}$$
In the right hand side $\tv_a(0)$ acts as $\td_a$, which gives
$$ J \tw  \otimes s \otimes v_h \otimes \td_a f,$$
while in the left hand side $\Phi_{ij} (\tv_a (0))$ acts as $\td_a + E_{sp} (0) \d_s \td_a x_p$. The left hand side then becomes
$$ J \tw  \otimes s \otimes v_h \otimes \td_a f 
+ (\td_a  J) \tw  \otimes s \otimes v_h \otimes  f 
+  E_{sp}(0) J \tw  \otimes s \otimes v_h \otimes (\d_s \td_a x_p) f .$$
In order to evaluate the action of $\td_a  J$ we note that
$(\td_a  J) J^{-1}$ belongs to the Lie algebra $\gl$, and we get
$$\td_a  J = ((\td_a  J)  J^{-1})  J
= (\td_a \d_s \tx_k) (\td_k x_p) E_{sp}(0) J .$$
Since $ (\td_a \d_s \tx_k) (\td_k x_p) + \d_s \td_a x_p = 0$, we establish (\vaz).

 Having established the homomorphism $\Psi_{ij} : \;  \M_0^\prime \rightarrow \M_0^{\prime\prime}$ as both $\ss_0$- and $\O_X(U_j \cap U_i)$-modules, we note that
$\M^\prime$ is a generalized Verma module for the Lie algebra $\ss$, generated by
$\M_0^\prime$, 
$$\M^\prime = U(\ss_-) \otimes \M_0^\prime,$$
thus, $\Psi_{ij}$ extends uniquely to a homomorphism
$$\Psi_{ij} : \;  \M^\prime \rightarrow \M^{\prime\prime}$$
of $\ss$-modules. By Corollary \Mmap, this is a homomorphism of modules for the vertex
algebra $\V(U_j \cap U_i)$. One can immediately see that $\Psi_{ji} \circ \Psi_{ij}$ is the identity map on 
\break
$\M(U_j \cap U_i)$, so $\Psi_{ij}$ is in fact an 
isomorphism of modules. 

 The $\V(U_i)$-module $\M(U_i)$ has a unique maximal submodule that trivially intersects with  $\M_0(U_i)$. The quotient module $\L(U_i)$ can be written as a tensor
product
$$\L(U_i) = \VHei \otimes \Lgl (W) \otimes \Laff(S) \otimes \LVir(h) \otimes \O_X(U_i),$$
where $\Lgl(W)$, $\Laff(S)$ and $\LVir(h)$ are the simple quotients of the 
corresponding 
\break
$\hgl$-, $\hg$- and Virasoro modules. It is clear that taking this quotient 
is compatible with the coordinate change map $\Psi_{ij}$, and we obtain a sheaf
$\L$ of modules for $\V$.

We established the following

{\bf Theorem \sheafL.} Let $\Lgl$ be an irreducible highest weight module for the Lie algebra
$\hgl$ at level $1$, such that its $\gl$-submodule $W$ generated the the highest weight 
vector of $\Lgl$ is a finite-dimensional rational $GL_N$-module. Let $\Laff$ be an 
irreducible highest weight module for $\hg$ at level $c \neq 0, -h^\vee$, and let $\LVir$ 
be an irreducible highest weight module for the Virasoro algebra with central charge
$ - \cvir$. There is a sheaf $\L$ of modules for the sheaf $\V$ of vertex algebras, where for an open set $U_i$ with a system of local parameters, 
the module $\L(U_i)$ is defined as
$$\L(U_i) = \VHei \otimes \Lgl  \otimes \Laff \otimes \LVir \otimes \O_X(U_i),$$ 
and the coordinate change map $\Psi_{ij}$ is defined on the top graded component
by (\mtop) and extended to $\L(U_j \cap U_i)$ as a homomorphism of 
$\ss$-modules.

\

{\bf 4. Representations of the sheaf $\G$ of Lie algebras.}

In this section we are going to show that the sheaf $\V$ of vertex algebras admits an 
action of the sheaf $\G$ of Lie algebras. As an immediate consequence we get 
representations of $\G$ on the sheaves of modules $\M$ and $\L$.

 Let $U \subset U_i$. For $f \in \O_X(U)$ we set the following formal  generating series
which coefficients span $\G(U)$:
$$\eqalign{
&g(f,z) = \sum_{j\in\Z} t^j f \otimes g z^{-j-1}, \;\; g \in \dg, \cr
&k_0(f,z) =  \sum_{j\in\Z} t^j f dt  z^{-j-1}, \cr
&k_a(f,z) =  \sum_{j\in\Z} t^j f dx_a  z^{-j-1}, \cr
&d_a(f,z) =  \sum_{j\in\Z} t^j f \d_a  z^{-j-1}, \;\; a =1, \ldots, N, \cr
&d_0(f,z) = - \sum_{j\in\Z} t^j f {\d \over \d t} z^{-j-1}. \cr}$$
The negative sign in the last formula is chosen to conform with the Virasoro algebra 
conventions.

{\bf Theorem \main.} 
Let  $V=\VHei \otimes \Vgl \otimes \Vaff \otimes \VVir$ be a tensor product of vertex
operator algebras, where $\Vgl$ is the universal enveloping vertex operator algebra for $\hgl$ at level $1$ and rank $-2N$, 
$\Vaff$ be the universal enveloping algebra for $\hg$ at
a non-zero, non-critical level $c$, and $\VVir$ be the universal enveloping Virasoro vertex algebra of rank $ - \cvir$, so that the total rank of $V$ is zero. Let $\V$ be the corresponding sheaf of vertex algebras on $X$. There is a representation $\rho$ of the sheaf of Lie algebras 
$$\G =  \left( \Map (\hX, \g) \oplus \bom^1_{\hX} \right) \rtimes \Vect(\hX)$$ 
on the sheaf of vertex algebras $\V$, given locally by the correspondence:
$$\rho(g(f,z)) = Y(g(-1) f,z), \eqno{(\Yg)}$$
$$\rho(k_0(f,z)) = c Y(f,z), \eqno{(\Yko)}$$
$$\rho(k_a(f,z)) = c Y(u_a(-1) f,z), \eqno{(\Yka)}$$
$$\rho(d_a(f,z)) = Y(v_a(-1) f,z) + \sum_{p=1}^N Y(E_{pa}(-1) \d_p f,z), \eqno{(\Yda)}$$
$$\rho(d_0(f,z)) = Y(\om_{(-1)} f,z) + \sum_{s,k=1}^N Y(u_k(-1) E_{sk}(-1) \d_s f,z)
- \sum_{p=1}^N Y(u_p(-2) \d_p f,z) . \eqno{(\Ydo)}$$

{\bf Proof.} We need to prove that everything is well-defined and that the Lie brackets
of the vertex operators in the right hand sides of (\Yg)-(\Ydo) match the Lie brackets
of the left hand sides.

 Note that the relation
$${d \over d z} Y(f,z) = \sum_{p=1}^N Y(u_p(-1) \d_p f, z)$$
ensures that the elements of $d (\C[t,t^{-1}] \O_X (U_i))$ act  trivially.

 We need to show that both sides of (\Yg)-(\Ydo) transform in a compatible way
under the coordinate changes in $U_i \cap U_j$, i.e.,
$$\Phi \circ \rho = \rho \circ \Theta ,$$
where $\Theta$ is the coordinate transformation
$\Theta : \; \G(U_j \cap U_i) \rightarrow \G(U_i \cap U_j)$.
For (\Yg) and (\Yko) this holds trivially. Let us verify this for (\Yka)-(\Ydo).
$$\Phi (\rho( \sum_{j\in\Z} t^j f d \tx_a z^{-j-1} )) = c \Phi (Y(\tu_a(-1) f,z))
= c Y((u_p(-1) \d_p \tx_a)_{(-1)} f, z) $$
$$= \rho( \sum_{j\in\Z} t^j f (\d_p \tx_a) d x_p z^{-j-1})
 = \rho (\Theta( \sum_{j\in\Z} t^j f d \tx_a z^{-j-1} )). $$
 Since $\Theta( \sum_{j\in\Z} t^j f \td_a z^{-j-1} ) =  
\sum_{j\in\Z} t^j (\td_a x_s) f \d_s z^{-j-1}$, the verification of compatibility for (\Yda)
amounts to checking the equality
$$\Phi \big(\tv_a (-1) f + \tE_{qa} (-1) \td_q f \big)
= v_s(-1) (\td_a x_s) f + E_{ks} (-1) \d_k ((\td_a x_s)f) . \eqno{(\Phida)}$$
Let us prove this equality:
$$\Phi  \big(\tv_a (-1) f + \tE_{qa} (-1) \td_q f \big)$$
$$= \big(v_p(-1) (\td_a x_p) + E_{ks}(-1) \d_k \td_a x_s\big)_{(-1)} f
+\big(E_{ks} (-1)(\d_k \tx_q) (\td_a x_s) +u_s(-1) \td_a \d_s \tx_q\big)_{(-1)}\td_q f$$
$$= (\td_a x_p)_{(-2)} v_p(0) f + v_p(-1) (\td_a x_p) f
+ E_{ks}(-1) (\d_k \td_a x_s) f + E_{ks} (-1)(\d_k \tx_q) (\td_a x_s) (\td_q f)$$
$$+ u_s (-1) ( \td_a \d_s \tx_q)  (\td_q f)$$
$$= u_s(-1) (\d_s \td_a x_p) (\d_p f) +  u_s (-1) ( \td_a \d_s \tx_q)  (\td_q f)
+ v_p(-1) (\td_a x_p) f$$
$$+ E_{ks}(-1) (\d_k \td_a x_s) f
+ E_{ks} (-1) (\td_a x_s)(\d_k  f)$$
$$=  v_p(-1) (\td_a x_p) f +  E_{ks} (-1) \d_k ((\td_a x_s)f) .$$
Here we used (\dcom) in the last step.

Finally, for (\Ydo) we need to show that
$$\Phi  \big( \tom_{(-1)} f + \tu_k (-1) \tE_{sk} (-1) \td_s f - \tu_s (-2) \td_s f  \big)$$
$$= \om_{(-1)} f + u_b (-1) E_{ab} (-1) \d_a f - u_a (-2) \d_a f . \eqno{(\vYdo)}$$

{\bf Lemma \twom.} $\Phi_{ij} (\tom) = \om$.

{\bf Proof.} The Virasoro element in $\V(U)$ is the sum of the Virasoro elements for the
tensor factors:
$$\om = \omhei + \omg + \omgl + \omvir .$$
The coordinate change map $\Phi$ does not affect the components $\Vaff$ and
$\VVir$, so  
$$\Phi(\tomg) = \omg, \;\;  \Phi(\tomvir) = \omvir. \eqno{(\tomo)}$$

One can verify that
$$\Phi(\tomhei) = \omhei - (\d_k \tx_p)_{(-3)} (\td_p x_k)
-  E_{sm}(-1) (\d_s \tx_a)_{(-2)}(\td_a x_m) \eqno{(\tomh)}$$
and 
$$\Phi(\tomgl) = \omgl +  (\d_k \tx_p)_{(-3)} (\td_p x_k)
+ E_{sm}(-1) (\d_s \tx_a)_{(-2)}(\td_a x_m).  \eqno{(\tomn)}$$
Adding (\tomo), (\tomn) and (\tomh) together we get the claim of the Lemma.

\ifnum \version=\longversion
Let us check (\tomh). Recall that
$$\omhei = v_p(-1) u_p(-1) \vac = v_p(-1)  (x_p)_{(-2)} \vac .$$
Then
$$\Phi(\tomhei) = \Phi_{ij}(\tv_a (-1) (\tx_a)_{(-2)} \vac) $$
$$= (v_m(-1) (\td_a x_m))_{(-1)} (\tx_a)_{(-2)} \vac
+ (E_{sm}(-1) (\d_s \td_a x_m))_{(-1)}  (\tx_a)_{(-2)} \vac$$
$$= v_m(-1) u_r(-1)  (\td_a x_m) (\d_r \tx_a)
+  (\td_a x_m)_{(-2)} v_m (0) (\tx_a)_{(-2)} \vac$$
$$+  (\td_a x_m)_{(-3)} v_m (1) (\tx_a)_{(-2)} \vac
+ E_{sm}(-1) u_r (-1) (\d_s \td_a x_m) (\d_r \tx_a)$$
$$=\omhei +  (\td_a x_m)_{(-2)} (\d_m \tx_a)_{(-2)} \vac
+  (\td_a x_m)_{(-3)} (\d_m \tx_a)_{(-1)} \vac$$
$$-  E_{sm}(-1) u_r (-1)  (\td_a x_m) (\d_s \d_r \tx_a)$$
$$= \omhei  -  (\td_a x_m)_{(-1)} (\d_m \tx_a)_{(-3)} \vac
- E_{sm} (-1) (\d_s \tx_a)_{(-2)}  (\td_a x_m).$$

Next let us derive (\tomn). The Virasoro element in $\Vgl$ is
$$\omgl = {1\over 2(N+1)} \left( I(-1) I(-1) \vac
+  \sum_{a,b=1}^N  E_{ab} (-1) E_{ba}(-1) \vac \right) + {1\over 2} I(-2) \vac .$$
We have
$$\Phi(\tI(-1) \vac) = E_{sp}(-1) (\d_s \tx_a) (\td_a x_p) - (\td_a x_p)_{(-2)}(\d_p \tx_a)$$
$$= I(-1) - (\td_a x_p)_{(-2)}(\d_p \tx_a), $$
$$\Phi(\tI(-2) \vac) = I(-2) \vac -  (D \td_a x_p)_{(-2)}(\d_p \tx_a)
 -  (\td_a x_p)_{(-2)}(D \d_p \tx_a)$$
$$= I(-2) \vac -2  (\td_a x_p)_{(-3)}(\d_p \tx_a) 
-  (\td_a x_p)_{(-2)}(\d_p \tx_a)_{(-2)}\vac ,$$
and
$$\Phi(\tomgl) = {1 \over 2(N+1)} \bigg( 
\Phi(\tI(-1))_{(-1)} \Phi(\tI(-1) \vac) 
+ \Phi(\tE_{ab}(-1))_{(-1)} \Phi(\tE_{ba}(-1))\vac\bigg) $$
$$+ {1\over 2}  \Phi(\tI(-2) \vac)$$
$$= {1 \over 2(N+1)}
 \bigg( 
\left(I(-1)\vac - (\td_a x_p)_{(-2)}(\d_p \tx_a)\right)_{(-1)} 
\left(I(-1) \vac - (\td_b x_s)_{(-2)}(\d_s \tx_b)\right)$$
$$
\kern-10pt
+ \big( E_{sp}(-1)(\d_s \tx_a) (\td_b x_p) - (\td_b x_p)_{(-2)} (\d_p \tx_a) \big)_{(-1)}
\kern-1pt
\big( E_{km}(-1)(\d_k \tx_b) (\td_a x_m) - (\td_a x_q)_{(-2)}(\d_q \tx_b)  \big)
\bigg)$$
$$+ {1 \over 2} I(-2) - (\td_a x_p)_{(-3)} (\d_p \tx_a) 
- {1 \over 2} (\td_a x_p)_{(-2)} (\d_p \tx_a)_{(-2)} \vac$$
$$={1 \over 2(N+1)} \bigg( 
 I(-1) I(-1) \vac - 2 I(-1) (\td_a x_p)_{(-2)} (\d_p  \tx_a)$$
$$+ (\td_a x_p)_{(-2)} (\td_b x_s)_{(-2)} (\d_p \tx_a) (\d_s \tx_b)
+ E_{sp}(-1) E_{km}(-1) (\d_s \tx_a) (\td_b x_p)(\d_k \tx_b)(\td_a x_m)$$
$$+ ((\d_s \tx_a)(\td_b x_p))_{(-2)} E_{sp}(0) E_{km}(-1) (\d_k \tx_b) (\td_a x_m)$$
$$+ ((\d_s \tx_a)(\td_b x_p))_{(-3)} E_{sp}(1) E_{km}(-1) (\d_k \tx_b) (\td_a x_m)$$
$$- E_{km}(-1) (\td_b x_p)_{(-2)} (\d_p \tx_a)(\d_k \tx_b)(\td_a x_m)
-E_{sp}(-1)(\td_a x_q)_{(-2)}  (\td_b x_p)(\d_q \tx_b)(\d_s \tx_a) $$
$$+ (\td_b x_p)_{(-2)} (\td_a x_q)_{(-2)} (\d_p \tx_a) (\d_q \tx_b)
 \bigg)$$
$$+ {1 \over 2} I(-2) \vac - (\td_a x_p)_{(-3)} (\d_p \tx_a) - {1\over 2} (\td_a x_p)_{(-2)} (\d_p \tx_a)_{(-2)} \vac$$
$$={1 \over 2(N+1)} \bigg( 
 I(-1) I(-1) \vac - 2 I(-1) (\td_a x_p)_{(-2)} (\d_p \tx_a)$$
$$+ (\td_a x_p)_{(-2)} (\td_b x_s)_{(-2)} (\d_p \tx_a) (\d_s \tx_b)
+ E_{sp}(-1) E_{ps}(-1) \vac$$
$$+ E_{sm}(-1) (\d_s \tx_a)_{(-2)} (\td_b x_p) (\td_a x_m)(\d_p \tx_b)
+ E_{sm}(-1) (\td_b x_p)_{(-2)} (\d_s \tx_a) (\td_a x_m) (\d_p \tx_b)$$
$$- E_{kp}(-1) (\d_s \tx_a)_{(-2)} (\td_b x_p) (\td_a x_s)(\d_k \tx_b)
- E_{kp} (-1) (\td_b x_p)_{(-2)} (\d_s \tx_a) (\td_a x_s)(\d_k \tx_b)$$
$$-2 E_{kp}(-1) (\td_b x_p)_{(-2)} (\d_k \tx_b)
+ (\d_s \tx_a)_{(-3)} (\td_b x_p) (\d_p \tx_b)(\td_a x_s)$$
$$+  (\td_b x_p)_{(-2)} (\d_s \tx_a)_{(-2)}(\d_p \tx_b) (\td_a x_s)
+ (\td_b x_p)_{(-3)} (\d_s \tx_a) (\d_p \tx_b)(\td_a x_s)$$
$$+ (\td_b x_p)_{(-2)} (\td_a x_q)_{(-2)} (\d_p \tx_a)(\d_q \tx_b)
 \bigg) $$
$$+ {1 \over 2} I(-2) \vac - (\td_a x_p)_{(-3)} (\d_p \tx_a) - {1\over 2} (\td_a x_p)_{(-2)} (\d_p \tx_a)_{(-2)} \vac$$
$$ = \omgl + 
{1 \over 2(N+1)} \bigg( 
- 2 I(-1) (\td_a x_p)_{(-2)} (\d_p \tx_a)
+ N E_{sm}(-1) (\d_s \tx_a)_{(-2)} (\td_a x_m)$$
$$+ I(-1) (\td_b x_p)_{(-2)} (\d_p \tx_b)
- I(-1) (\d_s \tx_a)_{(-2)}(\td_a x_s)$$
$$-N E_{kp}(-1) (\td_b x_p)_{(-2)} (\d_k \tx_b)
- 2 E_{kp}(-1) (\td_b x_p)_{(-2)} (\d_k \tx_b)
+ N (\d_s \tx_a)_{(-3)} (\td_a x_s)$$
$$+ (\td_b x_p)_{(-2)} ((\d_s \tx_a)(\td_a x_s))_{(-2)} (\d_p \tx_b)
+ N (\td_b x_p)_{(-3)} (\d_p \tx_b)$$
$$- (\td_b x_p)_{(-2)} (\d_p \tx_a)_{(-2)} (\td_a x_q) (\d_q \tx_b)
 \bigg)
- (\td_a x_p)_{(-3)} (\d_p \tx_a) 
 - {1\over 2} (\td_a x_p)_{(-2)} (\d_p \tx_a)_{(-2)} \vac$$
$$ = \omgl 
+ E_{sm}(-1) (\d_s \tx_a)_{(-2)} (\td_a x_m)
 - {1\over 2} (\td_a x_p)_{(-2)} (\d_p \tx_a)_{(-2)} \vac$$
$$- (\td_a x_p)_{(-3)} (\d_p \tx_a)
 - {1\over 2} (\td_a x_p)_{(-2)} (\d_p \tx_a)_{(-2)} \vac$$
$$ = \omgl 
+ E_{sm}(-1) (\d_s \tx_a)_{(-2)} (\td_a x_m)
+ (\d_p \tx_a)_{(-3)} (\td_a x_p) .$$
\fi

Now let us establish (\vYdo):
$$\Phi \left( \tom_{(-1)} f + \tu_k (-1) \tE_{sk} (-1) \td_s f - \tu_s (-2) \td_s f \right)$$
$$= \om_{(-1)} f + (u_b(-1) (\d_b \tx_k))_{(-1)} 
\left( E_{ac}(-1)(\d_a \tx_s) (\td_k x_c) + u_a (-1) \td_k \d_a \tx_s \right)_{(-1)}
(\td_s f) $$
$$- (u_a(-1) (\d_a \tx_s) )_{(-2)} (\td_s f)$$
$$ = \om_{(-1)} f 
+ u_b(-1) E_{ac} (-1) (\d_b \tx_k) (\d_a \tx_s) (\td_k x_c) (\td_s f)
+ u_b(-1) u_a(-1) (\d_b \tx_k) (\td_k \d_a \tx_s) (\td_s f)$$
$$- u_a(-2) (\d_a \tx_s)(\td_s f)
- u_a (-1) u_b(-1) (\d_b \d_a \tx_s) (\td_s f)$$
$$=  \om_{(-1)} f 
+ u_b (-1) E_{ab}(-1) (\d_a f)
+ u_b(-1) u_a (-1) (\d_b \d_a \tx_s) (\td_s f)$$
$$- u_a (-1) u_b(-1) (\d_b \d_a \tx_s) (\td_s f)
- u_a (-2) (\d_a f)$$
$$=  \om_{(-1)} f 
+ u_b (-1) E_{ab}(-1) (\d_a f)
- u_a (-2) (\d_a f). $$

 To complete the proof of Theorem \main, we need to show that locally (\Yg)-(\Ydo)
define a representation of the Lie algebra $\G(U_i)$.

 The Lie bracket in $\G(U_i)$ may be encoded using the commutators of the 
generating series:
$$[g_1 (f, z_1), g_2 (h, z_2)] = 
[g_1, g_2] (f h, z_2) \dz$$
$$+ (g_1 | g_2) k_0 (f h, z_2) \ddz
+ (g_1 | g_2) \sum_{p=1}^N k_p (h \d_p f, z_2) \dz, \eqno{(\Cgg)}$$
$$[d_a (f, z_1), g(h, z_2)] =  g(f \d_a h, z_2) \dz, \eqno{(\Cdag)}$$
$$[d_a (f, z_1), k_0 (h, z_2)] =  k_0 (f \d_a h, z_2) \dz, \eqno{(\Cdako)}$$ 
$$[d_a (f, z_1), k_b (h, z_2)] =  k_b (f \d_a h, z_2) \dz$$
$$+ \delta_{ab} k_0 (f h, z_2) \ddz
+ \delta_{ab} \sum_{p=1}^N k_p (h \d_p f, z_2) \dz, \eqno{(\Cdakb)}$$
$$[d_a (f, z_1), d_b(h, z_2)] = 
\big( d_b (f \d_a h, z_2) - d_a (h \d_b f, z_2) \big) \dz,  \eqno{(\Cdadb)}$$
$$[k_i (f, z_1), k_j (h, z_2)] = 0, \eqno{(\Ckikj)}$$
$$[g (f, z_1), k_i (h, z_2)] = 0, \eqno{(\Cgki)}$$
$$[d_0 (f, z_1), g(h, z_2)] =
\left( {\d \over \d z_2} g(f h, z_2) \right) \dz
+ g(f h, z_2) \ddz, \eqno{(\Cdog)}$$
$$[d_0 (f, z_1), k_0(h, z_2)] =
\sum_{p=1}^N k_p(f \d_p h, z_2) \dz,  \eqno{(\Cdoko)}$$
$$[d_0 (f, z_1), k_a(h, z_2)] =
\left( {\d \over \d z_2} k_a(f h, z_2) \right) \dz
+ k_a(f h, z_2) \ddz,  \eqno{(\Cdoka)}$$
$$[d_0 (f, z_1), d_a(h, z_2)] =
\left( {\d \over \d z_2} d_a(f h, z_2) \right) \dz$$
$$+ d_a(f h, z_2) \ddz
- d_0( h \d_a f, z_2) \dz, \eqno{(\Cdoda)}$$
$$[d_0 (f, z_1), d_0(h, z_2)] =
\left( {\d \over \d z_2} d_0(f h, z_2) \right) \dz
+ 2 d_0(f h, z_2) \ddz, \eqno{(\Cdodo)}$$
where $g, g_1, g_2 \in \dg$, $f, h \in \O_X(U_i)$, $a,b = 1, \ldots, N$, 
$i,j = 0, \ldots, N$.

 We will use the commutator formula (\comm) in order to prove that $\rho$ preserves these relations. For (\Cgg) we need to verify in $\V(U_i)$ the following relations
for $n$-th products:
$$ (g_1 (-1) f )_{(0)} (g_2 (-1) h) = 
[g_1, g_2](-1) f h 
+ (g_1 | g_2) c u_p (-1) h \d_p f ,$$
$$(g_1 (-1) f )_{(1)} (g_2 (-1) h) = 
c (g_1 | g_2) f h ,$$
$$(g_1 (-1) f )_{(n)} (g_2 (-1) h) = 0 {\;\; {\rm for} \;\;} n > 1.$$
Let us check these equalities:
$$ (g_1 (-1) f )_{(0)} (g_2 (-1) h) = 
f_{(-1)} g_1 (0) g_2 (-1) h
+ f_{(-2)} g_1 (1) g_2 (-1) h$$
$$= [g_1, g_2] (-1) f h + c  (g_1 | g_2)  u_p (-1) h \d_p f, $$
$$(g_1 (-1) f )_{(1)} (g_2 (-1) h) =
f_{(-1)} g_1 (1) g_2 (-1) h
= c (g_1 | g_2) f h .$$
The last relation holds trivially since the degree of the left hand side becomes negative.

 We verify (\Cdag) in an analogous way:
$$\big( v_a (-1) f + E_{pa} (-1) \d_p f \big)_{(0)} g(-1) h
= f_{(-1)} v_a (0) g(-1) h 
= g(-1) f \d_a h .$$
It is easy to see that
$$\big( v_a (-1) f + E_{pa} (-1) \d_p f \big)_{(n)} g(-1) h = 0
{\;\; \rm for \;\;} n > 0.$$

 For the remaining $n$-th products we will verify only those that have a non-negative degree.   For (\Cdako) we have
$$\big( v_a (-1) f + E_{pa} (-1) \d_p f \big)_{(0)} h
= f_{(-1)} v_a (0)  h  = f \d_a h .$$
To prove (\Cdakb), we compute
$$\big( v_a (-1) f + E_{pa} (-1) \d_p f \big)_{(0)} u_b(-1) h
=  f_{(-1)} v_a (0) u_b(-1) h
+  f_{(-2)} v_a (1) u_b(-1) h$$
$$= u_b (-1)  f \d_a h 
+ \delta_{ab} u_p (-1) h \d_p f ,$$
and
$$\big( v_a (-1) f + E_{pa} (-1) \d_p f \big)_{(1)} u_b(-1) h
=  f_{(-1)} v_a (1) u_b(-1) h 
= \delta_{ab} f h .$$

 Let us now verify (\Cdadb):
$$\big( v_a (-1) f + E_{pa} (-1) \d_p f \big)_{(0)}
\big( v_b (-1) h + E_{kb} (-1) \d_k h \big)$$
$$ = f_{(-1)} v_a (0) \big( v_b (-1) h + E_{kb} (-1) \d_k h \big)
+ v_a (-1) f_{(0)} v_b (-1) h $$
$$+ ( \d_p f)_{(-1)} E_{pa} (0)  E_{kb} (-1) \d_k h 
+ ( \d_p f)_{(-2)} E_{pa} (1)  E_{kb} (-1) \d_k h 
+  E_{pa} (-1)  ( \d_p f)_{(0)}  v_b (-1) h $$
$$=  f_{(-1)} v_b (-1) \d_a h
+ E_{kb}(-1) f \d_a \d_k h
- v_a(-1) (\d_b f)_{(-1)} h$$
$$+  ( \d_p f)_{(-1)} E_{pb} (-1) \d_ a h
-  ( \d_b f)_{(-1)} E_{ka} (-1) \d_ k h
+  ( \d_b f)_{(-2)} \d_a h
- E_{pa} (-1) (\d_b \d_p f ) h$$
$$= v_b (-1) f \d_ a h
-  ( \d_b f)_{(-2)} \d_a h
+  E_{kb}(-1) f \d_a \d_k h
+ E_{kb}(-1) (\d_k f) (\d_a  h)$$
$$- v_a (-1) (\d_b f) h 
- E_{ka} (-1) (\d_b f) (\d_k h)
- E_{ka} (-1) (\d_b \d_k f) h
+ (\d_b f)_{(-2)} (\d_a h)$$
$$= v_b (-1)  f \d_ a h +  E_{kb}(-1) \d_k ( f \d_ a h)
- v_a (-1) (\d_b f) h - E_{ka} (-1) \d_k ( (\d_b f) h) ,$$
and for $n=1$:
$$\big( v_a (-1) f + E_{pa} (-1) \d_p f \big)_{(1)}
\big( v_b (-1) h + E_{kb} (-1) \d_k h \big)$$
$$=  f_{(0)} v_a (0) \big( v_b (-1) h + E_{kb} (-1) \d_k h \big)
+ ( \d_p f)_{(-1)} E_{pa} (1)  E_{kb} (-1) \d_k h$$
$$=  f_{(0)} v_b (-1) \d_a h + (\d_b f)_{(-1)} \d_a h
= - (\d_b f)_{(-1)} \d_a h + (\d_b f)_{(-1)} \d_a h = 0.$$

 Relations (\Ckikj) and (\Cgki) hold trivially.

 For the commutators involving $d_0 (f,z)$, we will be using the 
properties (\omd) and (\omdeg) of the Virasoro element.
 We will also need the following commutator relations (see [2]):
$$[\om_{(n)}, f_{(m)}] = - (n+m) f_{(n+m-1)},$$
$$[\om_{(n)}, g(m)] = -m g(n+m-1),$$
$$[\om_{(n)}, u_a(m)] = -m u_a(n+m-1),$$
$$[\om_{(n)}, v_a(m)] = -m v_a(n+m-1),$$
$$[\om_{(n)}, E_{ab}(m)] = -m E_{ab}(n+m-1) - \delta_{ab} \delta_{n+m-1,0} 
{n(n-1) \over 2} \Id ,$$
$$[\om_{(n)},\om_{(m)}] = (n-m) \om_{(n+m-1)} .$$
For (\Cdog) we get
$$\big( \om_{(-1)} f + u_k(-1) E_{sk} (-1) \d_s f - u_k (-2) \d_k f \big)_{(0)}
g(-1) h$$
$$= f_{(-1)} \om_{(0)} g(-1) h
+  f_{(-2)} \om_{(1)} g(-1) h$$
$$= f D(g(-1) h) + g(-1)  (D f)_{(-1)} h
= D(g(-1) f h )$$
and
$$\big( \om_{(-1)} f + u_k(-1) E_{sk} (-1) \d_s f - u_k (-2) \d_k f \big)_{(1)}
g(-1) h$$
$$= f_{(-1)} \om_{(1)} g(-1) h
= g(-1) f h .$$
For (\Cdoko) we have
$$\big( \om_{(-1)} f + u_k(-1) E_{sk} (-1) \d_s f - u_k (-2) \d_k f \big)_{(0)} h$$
$$= f_{(-1)} \om_{(0)} h +  f_{(-2)} \om_{(1)} h
= f_{(-1)} Dh = u_p (-1) f \d_p h,$$
and
$$\big( \om_{(-1)} f + u_k(-1) E_{sk} (-1) \d_s f - u_k (-2) \d_k f \big)_{(1)} h$$
$$= f_{(0)} \om_{(0)} h +  f_{(-1)} \om_{(1)} h = 0.$$

To verify that (\Cdoka) holds for $\rho$, we calculate:
$$\big( \om_{(-1)} f + u_k(-1) E_{sk} (-1) \d_s f - u_k (-2) \d_k f \big)_{(0)} u_a(-1) h$$
$$= f_{(-1)} \om_{(0)} u_a(-1) h 
+  f_{(-2)} \om_{(1)} u_a(-1)h
+  f_{(-3)} \om_{(2)} u_a(-1)h$$
$$= f_{(-1)} D(u_a(-1) h)
+ (Df)_{(-1)} u_a(-1) h
+  f_{(-3)}  u_a(0)h
= D(u_a(-1) f h),$$
for $n=1$:
$$\big( \om_{(-1)} f + u_k(-1) E_{sk} (-1) \d_s f - u_k (-2) \d_k f \big)_{(1)} u_a(-1) h$$
$$= f_{(0)} \om_{(0)} u_a(-1) h 
+  f_{(-1)} \om_{(1)} u_a(-1)h
+  f_{(-2)} \om_{(2)} u_a(-1)h$$
$$= f_{(0)} D(u_a(-1) h)
+ f_{(-1)} u_a(-1) h
+  f_{(-2)}  u_a(0)h
= u_a(-1) f h,$$
and for $n=2$:
$$\big( \om_{(-1)} f + u_k(-1) E_{sk} (-1) \d_s f - u_k (-2) \d_k f \big)_{(2)} u_a(-1) h$$
$$=  f_{(-1)} \om_{(2)} u_a(-1)h =  f_{(-1)}  u_a(0) h = 0.$$

\ifnum \version=\longversion
To establish (\Cdoda) we need to compute the corresponding $n=0, 1$ and $2$ products:
$$\big( \om_{(-1)} f + u_k(-1) E_{sk} (-1) \d_s f - u_k(-2) \d_k f\big)_{(0)}
\big(v_a(-1) h + E_{pa} (-1) \d_p h \big) $$
$$= f_{(-1)} \om_{(0)} \big(v_a(-1) h + E_{pa} (-1) \d_p h \big)
+ f_{(-2)} \om_{(1)} \big(v_a(-1) h + E_{pa} (-1) \d_p h \big)$$
$$+ f_{(-3)} \om_{(2)} \big(v_a(-1) h + E_{pa} (-1) \d_p h \big)
+ \om_{(-1)} f_{(0)} v_a (-1) h$$
$$+ (E_{sk}(-1) \d_s f)_{(-2)} u_k(1)  v_a(-1) h
+ u_k (-1) (E_{sk}(-1) \d_s f)_{(0)} \big(v_a(-1) h + E_{pa} (-1) \d_p h \big)$$
$$+  u_k (-2) (E_{sk}(-1) \d_s f)_{(1)} E_{pa} (-1) \d_p h 
+ 2 (\d_k f)_{(-3)} u_k(1) v_a(-1) h
- u_k(-2)  (\d_k f)_{(0)} v_a(-1) h$$
$$= f_{(-1)} D \big(v_a(-1) h + E_{pa} (-1) \d_p h \big)
+ (Df)_{(-1)}  (v_a(-1) h + E_{pa} (-1) \d_p h )$$
$$+  f_{(-3)} v_a(0) h -   f_{(-3)} \d_a h
 - \om_{(-1)} (\d_a f) h
+ \big( D(E_{pa}(-1) \d_p f) \big)_{(-1)} h$$
$$+ u_k(-1) E_{sk}(-1) (\d_s f)_{(0)} v_a(-1) h
+ u_k(-1) (\d_s f)_{(-1)} E_{sk}(0) E_{pa}(-1) \d_p h$$
$$+ u_k(-1) (\d_s f)_{(-2)} E_{sk}(1) E_{pa}(-1) \d_p h
+ u_k(-2) (\d_s f)_{(-1)} E_{sk}(1) E_{pa}(-1) \d_p h$$
$$+ 2 (\d_a f)_{(-3)} h
+ u_k (-2) (\d_a \d_k f) h$$
$$= D( f_{(-1)} v_a (-1) h) 
+ D( f_{(-1)}  E_{pa}(-1) \d_p h)
+ f_{(-3)} \d_a h 
-  f_{(-3)} \d_a h 
- \om_{(-1)} (\d_a f) h$$
$$+ \big( D(E_{pa}(-1) \d_p f) \big)_{(-1)} h
- u_k(-1) E_{sk} (-1) (\d_a \d_s f) h
+ u_k(-1) (\d_s f)_{(-1)} E_{sa} (-1) \d_k h$$
$$- u_k(-1) (\d_a  f)_{(-1)} E_{pk} (-1) \d_p h
+ u_k(-1) (\d_a f)_{(-2)} \d_k h$$
$$+ u_k (-2) (\d_a f) (\d_k h)
+ 2 (\d_a f)_{(-3)} h
+ u_k(-2) (\d_k \d_a f) h$$
$$= D(v_a (-1) fh)
- D( (\d_a f)_{(-2)} h)
+ D( E_{pa}(-1) f \d_p h)
- \om_{(-1)}  (\d_a f) h$$
$$+ \big( D( E_{pa} (-1) \d_p f) \big)_{(-1)} h
+ E_{pa} (-1) (\d_p f) D(h)
- u_k (-1) E_{sk} (-1) (\d_s \d_a f) h$$
$$- u_k (-1) E_{sk}(-1) (\d_a f) (\d_s h)
+ (\d_a f)_{(-2)} D h
+ u_k (-2) \d_k ((\d_a f) h)
+ 2  (\d_a f)_{(-3)} h$$
$$= D\big( v_a(-1) fh +E_{pa} (-1) \d_p (fh) \big)$$
$$- \big( \om_{(-1)} (\d_a f) h + u_k (-1) E_{sk} (-1) \d_s ((\d_a f) h)
- u_k (-2) \d_k ((\d_a f) h)\big).$$
For $n=1$:
$$\big( \om_{(-1)} f + u_k(-1) E_{sk} (-1) \d_s f - u_k(-2) \d_k f\big)_{(1)}
\big(v_a(-1) h + E_{pa} (-1) \d_p h \big) $$
$$= f_{(0)} \om_{(0)}  \big(v_a(-1) h + E_{pa} (-1) \d_p h \big)
+  f_{(-1)} \om_{(1)}  \big(v_a(-1) h + E_{pa} (-1) \d_p h \big)$$
$$+  f_{(-2)} \om_{(2)}  \big(v_a(-1) h + E_{pa} (-1) \d_p h \big)
+ (E_{sk} (-1) \d_s f)_{(-1)} u_k (1) v_a(-1) h$$
$$+ u_k (-1)  (E_{sk} (-1) \d_s f)_{(1)} \big(v_a(-1) h+ E_{pa} (-1) \d_p h \big)
+ 2 (\d_k f)_{(-2)} u_k (1) v_a (-1) h$$
$$=  f_{(0)} v_a (-2) h + f_{(0)} v_a (-1) Dh
+  f_{(-1)} v_a (-1) h + f_{(-1)} E_{pa} (-1) \d_p h $$
$$+ f_{(-2)} v_a (-0) h
 -  f_{(-2)} \d_a h
+ E_{sa} (-1) (\d_s f) h$$
$$+ u_k(-1) (\d_s f)_{(-1)} E_{sk} (1) E_{pa} (-1) \d_p h
+ 2 (\d_a f)_{(-2)} h$$
$$= - (\d_a f)_{(-2)} h
- (\d_a f)_{(-1)} D h
+ v_a (-1) fh 
- (\d_a f)_{(-2)} h
+ E_{pa} (-1) f \d_p h$$
$$+ f_{(-2)} \d_a h 
- f_{(-2)} \d_a h
+ E_{pa} (-1) (\d_p f) h
+ u_k (-1) (\d_a f) (\d_k h)
+ 2 (\d_a f)_{(-2)} h$$
$$ = v_a (-1) fh + E_{pa}(-1) \d_p (fh), $$
and for $n=2$:
$$\big( \om_{(-1)} f + u_k(-1) E_{sk} (-1) \d_s f - u_k(-2) \d_k f \big)_{(2)}
\big(v_a(-1) h + E_{pa} (-1) \d_p h \big) $$
$$= f_{(1)} \om_{(0)}  (v_a(-1) h + E_{pa} (-1) \d_p h )
+ f_{(0)} \om_{(1)}  (v_a(-1) h + E_{pa} (-1) \d_p h )$$
$$+ f_{(-1)} \om_{(2)}  \big(v_a(-1) h + E_{pa} (-1) \d_p h \big)
+  (E_{sk} (-1) \d_s f)_{(0)} u_k (1) v_a(-1) h$$
$$+ 2 (\d_k f)_{(-1)} u_k (1) v_a (-1) h$$
$$=  f_{(1)} v_a (-2) h + f_{(1)} v_a (-1) Dh
+  f_{(0)} v_a (-1) h
+  f_{(-1)} v_a (0) h
-  f_{(-1)} \d_a h
+ 2 (\d_a f) h$$
$$= - (\d_a f) h 
 - (\d_a f)_{(0)} D h 
- (\d_a f) h
+ f \d_a h - f \d_a h
+ 2 (\d_a f) h = 0,$$
which proves (\Cdoda).

 In the following computation, which establishes (\Cdodo), we will be using (\Da)
and the Borcherds identity (\qas).

 We begin with calculating $n=0$ product of the elements of the vertex algebra
corresponding to $d_0 (f,z)$ and $d_0(h,z)$:
\else
To establish (\Cdodo) one needs to verify $n=0,1,2,3$ products for the elements
of the vertex algebra corresponding to $d_0 (f,z)$ and $d_0(h,z)$. In particular 
for $n=0$ product we must prove the equality
$$\big(\om_{(-1)} f + u_k (-1) E_{sk}(-1) \d_s f - u_k (-2) \d_k f\big)_{(0)}
\big(\om_{(-1)} h + u_p (-1) E_{mp} (-1) \d_m h - u_p (-2) \d_p h \big)$$
$$ =D \big( \om_{(-1)} f h + u_p (-1) E_{mp} (-1) \d_m (f h) - u_p (-2) \d_p (fh) \big) .$$
We verify this equality below and leave analogous computations for $n=1,2,3$, as
well as verification of (\Cdoda) as an exercise for the reader.
\fi
$$\big(\om_{(-1)} f + u_k (-1) E_{sk}(-1) \d_s f - u_k (-2) \d_k f\big)_{(0)}
\big(\om_{(-1)} h + u_p (-1) E_{mp} (-1) \d_m h - u_p (-2) \d_p h \big)$$
$$= f_{(-1)} \om_{(0)}\big(\om_{(-1)} h + u_p (-1) E_{mp} (-1) \d_m h - u_p (-2) \d_p h \big)$$
$$+ f_{(-2)} \om_{(1)}\big(\om_{(-1)} h + u_p (-1) E_{mp} (-1) \d_m h - u_p (-2) \d_p h \big)$$
$$+ f_{(-3)} \om_{(2)}\big(\om_{(-1)} h + u_p (-1) E_{mp} (-1) \d_m h - u_p (-2) \d_p h \big)$$
$$+f_{(-4)} \om_{(3)}\big(\om_{(-1)} h + u_p (-1) E_{mp} (-1) \d_m h - u_p (-2) \d_p h \big)$$
$$+ \om_{(-1)} f_{(0)} \om_{(-1)} h
+ \om_{(-2)} f_{(1)} \om_{(-1)} h$$
$$+ (E_{sk}(-1) \d_s f)_{(-2)} u_k (1)  \om_{(-1)} h
+ (E_{sk}(-1) \d_s f)_{(-3)} u_k (2)  \om_{(-1)} h$$
$$+ u_k (-1)   (E_{sk}(-1) \d_s f)_{(0)} \big(\om_{(-1)} h + u_p (-1) E_{mp} (-1) \d_m h \big)$$
$$+ u_k (-2)   (E_{sk}(-1) \d_s f)_{(1)} \big(\om_{(-1)} h + u_p (-1) E_{mp} (-1) \d_m h\big)$$
$$+ u_k (-3)   (E_{sk}(-1) \d_s f)_{(2)} \big(\om_{(-1)} h + u_p (-1) E_{mp} (-1) \d_m h \big)$$
$$+ 3 (\d_k f)_{(-4)} u_k (2)  \om_{(-1)} h
+ 2 (\d_k f)_{(-3)} u_k (1)  \om_{(-1)} h$$
$$- u_k (-2)  (\d_k f)_{(0)} \om_{(-1)} h
- 2 u_k (-3)  (\d_k f)_{(1)} \om_{(-1)} h$$
$$=  f_{(-1)} D\big(\om_{(-1)} h + u_p (-1) E_{mp} (-1) \d_m h - u_p (-2) \d_p h\big)$$
$$+ 2 (Df)_{(-1)} \big(\om_{(-1)} h + u_p (-1) E_{mp} (-1) \d_m h - u_p (-2) \d_p h \big)$$
$$+ 3  f_{(-3)} \om_{(0)} h
+ f_{(-3)} u_p (0) E_{mp} (-1) \d_m h
+  f_{(-3)} u_p (-1)  \om_{(2)} E_{mp} (-1) \d_m h$$
$$- 2  f_{(-3)} u_p (-1)  \d_p h
+ 4 f_{(-4)} \om_{(1)} h
+  f_{(-4)} u_p (1) E_{mp} (-1) \d_m h
- 2  f_{(-4)} u_p (0) \d_p h$$
$$- \om_{(-1)}  f_{(-2)} h
+  (E_{sk}(-1) \d_s f)_{(-2)} u_k (-1) h
+ 2  (E_{sk}(-1) \d_s f)_{(-3)} u_k (0) h$$
$$+  u_k (-1)  (\d_s f)_{(-1)} E_{sk}(0)  \big(\om_{(-1)} h + u_p (-1) E_{mp} (-1) \d_m h \big)$$
$$+  u_k (-1)  (\d_s f)_{(-2)} E_{sk}(1)  \big(\om_{(-1)} h + u_p (-1) E_{mp} (-1) \d_m h \big)
+  u_k (-1)  (\d_s f)_{(-3)} E_{sk}(2)  \om_{(-1)} h $$
$$ +  u_k (-1)  E_{sk}(-1)  (\d_s f)_{(0)} \om_{(-1)} h 
+  u_k (-1)  E_{sk}(-2)  (\d_s f)_{(1)} \om_{(-1)} h $$
$$+  u_k (-2)  (\d_s f)_{(0)} E_{sk}(0)  \big(\om_{(-1)} h + u_p (-1) E_{mp} (-1) \d_m h\big)$$
$$+  u_k (-2)  (\d_s f)_{(-1)} E_{sk}(1)  \big(\om_{(-1)} h + u_p (-1) E_{mp} (-1) \d_m h \big)
+  u_k (-2)  (\d_s f)_{(-2)} E_{sk}(2)  \om_{(-1)} h $$
$$ +  u_k (-2)  E_{sk}(-1)  (\d_s f)_{(1)} \om_{(-1)} h $$
$$ +  u_k (-3)  (\d_s f)_{(1)} E_{sk}(0)  \big(\om_{(-1)} h + u_p (-1) E_{mp} (-1) \d_m h \big)$$
$$+  u_k (-3)  (\d_s f)_{(0)} E_{sk}(1)  \big(\om_{(-1)} h + u_p (-1) E_{mp} (-1) \d_m h \big)
+  u_k (-3)  (\d_s f)_{(-1)} E_{sk}(2)  \om_{(-1)} h $$
$$+ 6 (\d_k f)_{(-4)} u_k (0) h
+ 2 (\d_k f)_{(-3)} u_k (-1) h
+ u_k (-2)  (\d_k f)_{(-2)} h$$
$$= D  (f_{(-1)} \om_{(-1)} h + u_p (-1) E_{mp} (-1) f \d_m h - u_p (-2) f \d_p h)$$
$$+  f_{(-2)} \om_{(-1)} h
+ u_p (-1) E_{mp} (-1)  f_{(-2)} \d_m h
- u_p (-2)   f_{(-2)} \d_p h$$
$$+ 3  f_{(-3)} Dh
-  f_{(-3)} u_p (-1) \d_p h
- 2  f_{(-3)} Dh$$
$$- \om_{(-1)}  f_{(-2)} h
+ (D (E_{sk}(-1) \d_s f))_{(-1)} u_k (-1) h$$
$$+ u_k (-1) u_p (-1) E_{sp} (-1) (\d_s f) (\d_k h)
-  u_k (-1) u_p (-1) E_{mk} (-1) (\d_p f) (\d_m h)$$
$$+ u_k (-1) (\d_s f)_{(-2)} E_{sk}(-1) h
+ u_k (-1) (\d_s f)_{(-2)} u_s (-1) \d_k h
 + u_k (-1) (\d_k f)_{(-3)} h$$
$$- u_k (-1)E_{sk}(-1) (\d_s f)_{(-2)} h
+  u_k (-2) (\d_s f)_{(-1)} E_{sk}(-1) h
+ u_k (-2) u_s (-1) (\d_s f) (\d_k h) $$
$$+ u_k (-2)  (\d_k f)_{(-2)} h
+ u_k (-3)  (\d_k f) h
+ 2  u_k (-1)  (\d_k f)_{(-3)} h
+  u_k (-2)  (\d_k f)_{(-2)} h$$
$$ =  D  \big( \om_{(-1)} f h + u_p (-1) E_{mp} (-1) f \d_m h - u_p (-2) f \d_p h \big)$$
$$- 2 D (f_{(-3)} h) 
- 3 f_{(-4)} h$$
$$ + (D (E_{sk}(-1) \d_s f))_{(-1)} u_k (-1) h
+  u_p (-1) E_{sp} (-1)  ( \d_s f) Dh
+ (D u_k (-1))_{(-1)} E_{sk}(-1)  ( \d_s f) h$$
$$+ u_s (-1)  (\d_s f)_{(-2)} Dh
+ 3 u_k (-1)  (\d_k f)_{(-3)} h
+ 2  u_k (-2)  (\d_k f)_{(-2)} h
+  u_k (-3)  (\d_k f) h$$
$$ =  D  \big( \om_{(-1)} f h + u_p (-1) E_{mp} (-1) f \d_m h - u_p (-2) f \d_p h \big)$$
$$+ D( u_p (-1) E_{sp} (-1) (\d_s f) h )
- 3 (Df)_{(-3)} h
- (Df)_{(-2)} Dh
+ u_p (-1) (\d_p f)_{(-2)} Dh$$
$$+ 3 u_p (-1)  (\d_p f)_{(-3)} h
+ 2  u_p (-2) (\d_p f)_{(-2)} h
+ u_p (-3) (\d_p f ) h$$
$$ =  D \big( \om_{(-1)} f h + u_p (-1) E_{mp} (-1) \d_m (f h) - u_p (-2) f \d_p h \big)$$
$$ - 3 (u_p (-1) \d_p f)_{(-3)} h 
-  (u_p (-1) \d_p f)_{(-2)} D h 
+ u_p (-1) (\d_p f)_{(-2)} D h$$
$$+ 3  u_p (-1)  (\d_p f)_{(-3)} h
+ 2  u_p (-2)  (\d_p f)_{(-2)} h
+  u_p (-3)  (\d_p f) h$$
$$ = D \big( \om_{(-1)} f h + u_p (-1) E_{mp} (-1) \d_m (f h) - u_p (-2) f \d_p h \big)$$
$$-3 u_p (-3)  (\d_p f) h
-3 u_p (-2)  (\d_p f)_{(-2)} h
-3 u_p (-1)  (\d_p f)_{(-3)} h$$
$$- u_p (-2) (\d_p f) Dh
-  u_p (-1) (\d_p f)_{(-2)} D h 
+  u_p (-1) (\d_p f)_{(-2)} D h $$
$$ + 3  u_p (-1) (\d_p f)_{(-3)} h 
 + 2  u_p (-2) (\d_p f)_{(-2)} h 
+ u_p (-3) (\d_p f) h $$
$$ = D \big( \om_{(-1)} f h + u_p (-1) E_{mp} (-1) \d_m (f h) - u_p (-2) f \d_p h \big)
- D (u_p (-2) (\d_p f) h)$$
$$ = D \big( \om_{(-1)} f h + u_p (-1) E_{mp} (-1) \d_m (f h) - u_p (-2) \d_p (fh) \big) .$$

\ifnum \version=\longversion
Let us do the computations for $n=1$ product in (\Cdodo):
$$\big(\om_{(-1)} f + u_k (-1) E_{sk}(-1) \d_s f - u_k (-2) \d_k f \big)_{(1)}
\big(\om_{(-1)} h + u_p (-1) E_{mp} (-1) \d_m h - u_p (-2) \d_p h \big)$$
$$= f_{(0)} \om_{(0)}\big(\om_{(-1)} h + u_p (-1) E_{mp} (-1) \d_m h - u_p (-2) \d_p h \big)$$
$$+ f_{(-1)} \om_{(1)}\big(\om_{(-1)} h + u_p (-1) E_{mp} (-1) \d_m h - u_p (-2) \d_p h\big)$$
$$+ f_{(-2)} \om_{(2)}\big(\om_{(-1)} h + u_p (-1) E_{mp} (-1) \d_m h - u_p (-2) \d_p h \big)$$
$$+f_{(-3)} \om_{(3)}\big(\om_{(-1)} h + u_p (-1) E_{mp} (-1) \d_m h - u_p (-2) \d_p h \big)$$
$$+ \om_{(-1)} f_{(1)} \om_{(-1)} h
+ (E_{sk}(-1) \d_s f)_{(-1)} u_k (1)  \om_{(-1)} h
+ (E_{sk}(-1) \d_s f)_{(-2)} u_k (2)  \om_{(-1)} h$$
$$+ u_k (-1)   (E_{sk}(-1) \d_s f)_{(1)} \big(\om_{(-1)} h + u_p (-1) E_{mp} (-1) \d_m h\big)$$
$$+ u_k (-2)   (E_{sk}(-1) \d_s f)_{(2)} \big(\om_{(-1)} h + u_p (-1) E_{mp} (-1) \d_m h \big)$$
$$+ 3 (\d_k f)_{(-3)} u_k (2)  \om_{(-1)} h
+ 2 (\d_k f)_{(-2)} u_k (1)  \om_{(-1)} h
- u_k (-2)  (\d_k f)_{(1)} \om_{(-1)} h$$
$$ = f_{(0)} \om_{(-2)} h 
+ f_{(0)} \om_{(-1)} Dh 
+ 2 f_{(-1)} \om_{(-1)} h$$
$$ + 2  f_{(-1)} u_p (-1) E_{mp} (-1) \d_m h
- 2   f_{(-1)} u_p (-2) \d_p h
+ 3  f_{(-2)} \om_{(0)} h
+  f_{(-2)} u_p (0) E_{mp} (-1) \d_m h$$
$$+  f_{(-2)} u_p (-1) \om_{(2)} E_{mp} (-1) \d_m h
- 2  f_{(-2)} u_p (-1) \d_p h
+ 4  f_{(-3)} \om_{(1)} h
+  f_{(-3)} u_p (1) E_{mp} (-1) \d_m h$$
$$ - 2 f_{(-3)} u_p (0) \d_p h
+ (E_{sk}(-1) \d_s f)_{(-1)} u_k (-1)  h
+ 2 (E_{sk}(-1) \d_s f)_{(-2)} u_k (0)  h$$
$$+ u_k (-1) (\d_s f)_{(0)} E_{sk}(0)  (\om_{(-1)} h + u_p (-1) E_{mp} (-1) \d_m h)$$
$$+ u_k (-1) (\d_s f)_{(-1)} E_{sk}(1)  (\om_{(-1)} h + u_p (-1) E_{mp} (-1) \d_m h) $$
$$ + u_k (-1)  (\d_s f)_{(-2)}  E_{sk}(2)  \om_{(-1)} h
+ u_k (-1) E_{sk}(-1)  (\d_s f)_{(1)} \om_{(-1)} h$$
$$+ u_k (-2)  (\d_s f)_{(1)}  E_{sk}(0) \om_{(-1)} h
+ u_k (-2)  (\d_s f)_{(0)}  E_{sk}(1) \om_{(-1)} h$$
$$+ u_k (-2)  (\d_s f)_{(-1)}  E_{sk}(2) \om_{(-1)} h
+ 6  (\d_k f)_{(-3)} u_k (0) h
+ 2  (\d_k f)_{(-2)} u_k (-1) h$$
$$= -2 f_{(-3)} h -  f_{(-2)} D h + 2  \om_{(-1)} f h -4 f_{(-3)} h$$
$$+ 2 u_p (-1)  E_{mp} (-1) f \d_m h
- 2  u_p (-2)  f \d_p h
+ 3  f_{(-2)} D h$$
$$- u_p (-1) f_{(-2)} \d_p h
- 2 u_p (-1) f_{(-2)} \d_p h
+ E_{sk}(-1) u_k (-1) (\d_s f) h$$
$$+ u_k (-1) (\d_s f)_{(-1)} E_{sk}(-1)  h
+ u_k (-1) (\d_s f)_{(-1)} u_s (-1) \d_k h
+ u_k (-1)  (\d_k f)_{(-2)} h$$
$$ +  u_k (-2)  (\d_k f) h
+ 2 u_k (-1)  (\d_k f)_{(-2)} h$$
$$ =  2  \om_{(-1)} f h
+ 2 u_k (-1) E_{sk}(-1) \d_s (f h)
- 3 (Df)_{(-2)} h$$
$$- 2 u_k (-2) f \d_k h
+ u_k (-2) (\d_k f) h
+ 3 u_k (-1) (\d_k f)_{(-2)} h$$
$$ =  2  \om_{(-1)} f h
+ 2 u_k (-1) E_{sk}(-1) \d_s (f h)
- 3 u_k (-2) (\d_k f) h$$
$$-3 u_k (-1) (\d_k f)_{(-2)} h
-2 u_k (-2) f \d_k h
+  u_k (-2) (\d_k f) h
+ 3 u_k (-1) (\d_k f)_{(-2)} h$$
$$= 2  \big( \om_{(-1)} f h
+  u_k (-1) E_{sk}(-1) \d_s (f h)
-  u_k (-2) \d_k ( f h) \big) .$$

For $n=2$:
$$\big(\om_{(-1)} f + u_k (-1) E_{sk}(-1) \d_s f - u_k (-2) \d_k f \big)_{(2)}
\big(\om_{(-1)} h + u_p (-1) E_{mp} (-1) \d_m h - u_p (-2) \d_p h \big)$$
$$= f_{(1)} \om_{(0)}\big(\om_{(-1)} h + u_p (-1) E_{mp} (-1) \d_m h - u_p (-2) \d_p h\big)$$
$$+ f_{(0)} \om_{(1)}\big(\om_{(-1)} h + u_p (-1) E_{mp} (-1) \d_m h - u_p (-2) \d_p h\big)$$
$$+ f_{(-1)} \om_{(2)}\big(\om_{(-1)} h + u_p (-1) E_{mp} (-1) \d_m h - u_p (-2) \d_p h\big)$$
$$+f_{(-2)} \om_{(3)}\big(\om_{(-1)} h + u_p (-1) E_{mp} (-1) \d_m h - u_p (-2) \d_p h \big)$$
$$+ (E_{sk}(-1) \d_s f)_{(0)} u_k (1)  \om_{(-1)} h
+ (E_{sk}(-1) \d_s f)_{(-1)} u_k (2)  \om_{(-1)} h$$
$$+ u_k (-1)   (E_{sk}(-1) \d_s f)_{(2)} \big(\om_{(-1)} h + u_p (-1) E_{mp} (-1) \d_m h\big)$$
$$+ 2 (\d_k f)_{(-1)} u_k (1)  \om_{(-1)} h
+ 3 (\d_k f)_{(-2)} u_k (2)  \om_{(-1)} h$$
$$ = f_{(1)} \om_{(-2)} h
+ f_{(1)} \om_{(-1)} Dh
+ 2 f_{(0)} \om_{(-1)} h
+ 3  f_{(-1)} \om_{(0)} h$$
$$+  f_{(-1)} u_p (0) E_{mp} (-1) \d_m h
+  f_{(-1)} u_p (-1) \om_{(2)} E_{mp} (-1) \d_m h
-2  f_{(-1)} u_p (-1)  \d_p h$$
$$+ 4  f_{(-2)} \om_{(1)} h
+  f_{(-2)} u_p (1) E_{mp} (-1) \d_m h
-2  f_{(-2)} u_p (0)  \d_p h$$
$$ + (E_{sk}(-1) \d_s f)_{(0)} u_k (-1) h
 + 2 (E_{sk}(-1) \d_s f)_{(-1)} u_k (0) h$$
$$+ u_k (-1)  (\d_s f)_{(1)} E_{sk}(0) \big(\om_{(-1)} h + u_p (-1) E_{mp} (-1) \d_m h\big)$$
$$+ u_k (-1)  (\d_s f)_{(0)} E_{sk}(1) \big(\om_{(-1)} h + u_p (-1) E_{mp} (-1) \d_m h \big)$$
$$+ u_k (-1)  (\d_s f)_{(-1)} E_{sk}(2) \om_{(-1)} h
+ 2 (\d_k f)_{(-1)} u_k (-1)  h
+ 6 (\d_k f)_{(-2)} u_k (0) h$$ 
$$ = -  f_{(-2)} h - 2  f_{(-2)} h + 3  f_{(-1)} Dh
-  f_{(-1)} u_p (-1) \d_p h$$
$$ - 2  f_{(-1)} Dh
+ u_k (-1) (\d_k f)_{(-1)} h + 2 f_{(-2)} h = 0.$$

And finally for $n=3$:
$$\big(\om_{(-1)} f + u_k (-1) E_{sk}(-1) \d_s f - u_k (-2) \d_k f \big)_{(3)}
\big(\om_{(-1)} h + u_p (-1) E_{mp} (-1) \d_m h - u_p (-2) \d_p h \big)$$
$$= f_{(2)} \om_{(0)}\big(\om_{(-1)} h + u_p (-1) E_{mp} (-1) \d_m h - u_p (-2) \d_p h\big)$$
$$+ f_{(1)} \om_{(1)}\big(\om_{(-1)} h + u_p (-1) E_{mp} (-1) \d_m h - u_p (-2) \d_p h\big)$$
$$+ f_{(0)} \om_{(2)}\big(\om_{(-1)} h + u_p (-1) E_{mp} (-1) \d_m h - u_p (-2) \d_p h\big)$$
$$+f_{(-1)} \om_{(3)}\big(\om_{(-1)} h + u_p (-1) E_{mp} (-1) \d_m h - u_p (-2) \d_p h\big)$$
$$+ (E_{sk}(-1) \d_s f)_{(1)} u_k (1)  \om_{(-1)} h
+ (E_{sk}(-1) \d_s f)_{(0)} u_k (2)  \om_{(-1)} h$$
$$+ 2 (\d_k f)_{(0)} u_k (1)  \om_{(-1)} h
+ 3 (\d_k f)_{(-1)} u_k (2)  \om_{(-1)} h$$
$$=  f_{(2)} \om_{(-2)} h
+ f_{(2)} \om_{(-1)} Dh
+ 2 f_{(1)} \om_{(-1)} h
+ 3 f_{(0)} \om_{(0)} h
+ 4 f_{(-1)} \om_{(1)} h$$
$$+  (E_{sk}(-1) \d_s f)_{(1)} u_k (-1) h
+ 2  (E_{sk}(-1) \d_s f)_{(0)} u_k (0) h$$
$$+ 2 (\d_k f)_{(0)} u_k (-1)  h
+ 6 (\d_k f)_{(-1)} u_k (0) h$$
$$= f_{(0)} Dh + 3  f_{(0)} Dh 
+ (\d_s f)_{(0)} E_{sk}(0)  u_k (-1) h
+ (\d_s f)_{(-1)} E_{sk}(1)  u_k (-1) h = 0.$$
\fi
This completes the proof of Theorem \main.

\

Let us discuss some applications of our results to the chiral de Rham complex,
constructed by Malikov-Schechtman-Vaintrob [9, 8]. Here we specialize to
the case $\dg = (0)$. Under this assumption, the sheaf $\V$ has a local description
$$\V (U_i) = \VHei \otimes \Vgl \otimes \VVir \otimes \O_X (U_i),$$
where $\VVir$ is the universal enveloping vertex algebra for $\Vir$ of rank $0$.
The sheaf $\V$ is a module for the sheaf $\Vect (\hX)$ of Lie algebras.
The chiral de Rham complex is a sheaf $\R$ of vertex superalgebras on $X$  with
a local description
$$\R(U_i) = \VHei \otimes \VZN \otimes \O_X (U_i),$$
where $\VZN$ is a vertex superalgebra associated with the standard Euclidean 
lattice $\Z^N$. Malikov-Schechtman-Vaintrob use a fermionic realization of $\VZN$
with $\VZN$ being $\Z$-graded by fermionic degree,
$$\VZN = \mathop\oplus\limits_{k=-\infty}^\infty \VZN^{(k)}.$$
Each component $\VZN^{(k)}$ is an irreducible highest weight
module for $\hgl$ at level $1$. This induces a $\Z$-grading on the chiral de Rham 
complex, $\R = \mathop\oplus\limits_{k=-\infty}^\infty \R^{(k)}$. The chiral differential is a map
$$d: \;\; \R^{(k)} \rightarrow \R^{(k+1)}$$
(see [9] for details).

{\bf Theorem \chir.} (i) Each component $\R^{(k)}$ of the chiral de Rham complex
is a module for the sheaf $\Vect (\hX)$ of Lie algebras of vector fields.

(ii) The differential $d: \;\; \R^{(k)} \rightarrow \R^{(k+1)}$ is a homomorphism
of modules for the  sheaf $\Vect (\hX)$.

{\bf Proof.} Applying Theorem \main \ with  $\dg = (0)$, $\Lgl = \VZN^{(k)}$ and
$\LVir$ being trivial 1-dimensional module for the Virasoro algebra, we obtain that
$\R^{(k)}$ is a module for the sheaf of Lie algebras 
$\bom^1_{\hX} \rtimes  \Vect (\hX)$, from which the first claim follows.

The proof of the second claim is completely analogous to Theorem 10.3 in [3],
where the case of a torus is considered.

{\bf Remark.} Note that although  $d: \;\; \R^{(k)} \rightarrow \R^{(k+1)}$ is a homomorphism of  $\Vect (\hX)$-modules, it does not commute with the action of  
$\bom^1_{\hX}$.

\

{\bf 5. Lie algebras associated with the field of rational functions on $X$.}

 In conclusion of the paper we outline a version of our construction in the setting 
of rational functions. Let $\C(X)$ be the field of rational functions on $X$ and let
$R = \C[t,t^{-1}] \otimes \C(X)$. Consider the corresponding $R$-modules $\Omega^1_R$ of 1-forms on $\hX$, and $\Vect_R$ of vector fields on $\hX$.
Formally these objects may be defined as direct limits over non-empty open subsets
$U \subset X$:
$$\Omega^1_R = \dirlim \Omega^1_{\hX} (U), 
\quad\quad
\Vect_R = \dirlim \Vect_{\hX} (U).$$
 We consider the Lie algebra
$$\G_R = R \otimes \dg \oplus \left( \Omega^1_R / d R \right) \oplus \Vect_R,$$
where the definition of the Lie bracket is completely analogous to one given in section 1.
Alternatively, $\G_R$ may be defined as a direct limit
$$\G_R = \dirlim \G(U) .$$
The field $\C(X)$ is a finite extension
of a purely transcendental extension of $\C$, 
$\C(X) = \C(x_1, \ldots, x_N; y_1, \ldots, y_s)$, where
$\{ x_1, \ldots, x_N \}$ are algebraically independent and $y_i$'s are algebraic over
$\C(x_1, \ldots, x_N)$. Note that $\Omega^1_R$ is a free $R$-module with the free generators $dt, dx_1, \ldots, dx_N$, and $\Vect_R$ is a free $R$-module with the
generators ${\d \over \d t}, {\d \over \d {x_1}}, \ldots, {\d \over \d {x_N}}$.

 We can use $\{ x_1, \ldots, x_N \}$ in place of the local coordinates to define the 
vertex algebra
$$\V_R = \VHei \otimes \Vgl \otimes \Vaff \otimes \VVir \otimes \C(X) .$$
Our construction yields the following result:

{\bf Theorem \ratl.}
Let $\V_R = \VHei \otimes \Vgl \otimes \Vaff \otimes \VVir \otimes \C(X) $ be the vertex
algebra with the central charges of the tensor factors the same as in Theorem \main.

(i) $\V_R$ is a module for the Lie algebra $\G_R$, where the action is given by (\Yg)-(\Ydo).

(ii) Let $\Lgl$, $\Laff$, $\LVir$ be irreducible modules for the vertex algebras
$\Vgl$, $\Vaff$, $\VVir$ respectively. Then
$$\L_R = \VHei \otimes \Lgl \otimes \Laff \otimes \LVir \otimes \C(X)$$
is an irreducible module for the Lie algebra $\G_R$.

{\bf Proof.} The only claim that requires a proof here is the irreducibility of $\L_R$.
Everything else follows immediately from Theorems \main \ and \chir \ by passing to
the direct limit. 

To show that $\L_R$ is irreducible as a $\G_R$-module, we note that the fields in
(\Yg)-(\Ydo) that define the action of $\G_R$ on $\V_R$ generate the vertex algebra
$\V_R$. Since $\Lgl$, $\Laff$, $\LVir$ are irreducible modules for the vertex algebras
$\Vgl$, $\Vaff$, $\VVir$ respectively, and the vertex algebra
$ \VHei \otimes \C(X)$ is simple, we conclude that $\L_R$ is an irreducible $\G_R$-module.

\

{\bf References:}

\noindent
[1] S.~Berman, Y.~Billig, 
{Irreducible representations for toroidal Lie algebras.} 
\break
J.Algebra {\bf 221} (1999), 188-231.

\noindent
[2] Y.~Billig, 
{A category of modules for the full toroidal Lie algebra.}
 Int. Math. Res. Not., (2006), Art. ID 68395, 46 pp.

\noindent
[3] Y.~Billig, V.~Futorny, 
{Representations of Lie algebra of vector fields on a torus
and chiral de Rham complex.} arXiv:1108.6092 [math.RT].

\noindent
[4] S.~Eswara Rao, R.V.~Moody, 
{Vertex representations for $n$-toroidal
Lie algebras and a generalization of the Virasoro algebra.}
Commun.Math.Phys. {\bf 159} (1994), 239-264.

\noindent
[5] V.~Kac, 
{Vertex Algebras for Beginners.}
University Lecture Series, {\bf 10},
Amer. Math. Soc, Providence,
2nd Edition, 1998.

\noindent
[6] T.~A.~Larsson,
{Lowest-energy representations of non-centrally extended diffeomorphism 
algebras.}
Commun.Math.Phys. {\bf 201} (1999), 461-470.

\noindent
[7] H.~Li, 
{Local systems of vertex operators, vertex superalgebras and modules.}
 J.Pure Appl.Algebra {\bf 109} (1996), 143-195.

\noindent
[8] F.~Malikov, V.~Schechtman,
{Chiral de Rham complex. II.}
Differential topology, infinite-dimensional Lie algebras, and applications, 149-188, Amer. Math. Soc. Transl. Ser. 2, 194, Amer. Math. Soc., Providence, RI, 1999. 

\noindent
[9] F.~Malikov, V.~Schechtman, A.~Vaintrob,
{Chiral de Rham complex.}
Commun. Math. Phys., {\bf 204} (1999), 439-473.

\noindent
[10] R.V.~Moody, S.E.~Rao, T.~Yokonuma,  
{Toroidal Lie algebras and vertex representations.} 
Geom.Ded. {\bf 35} (1990), 283-307.

\noindent
[11] D.~Mumford, 
{The Red Book of Varieties and Schemes.}
Lecture Notes in Mathematics, {\bf 1358}, Springer-Verlag: Berlin/New York, 1988.

\end